\documentclass[11pt]{amsart}
\usepackage{a4wide}
\usepackage{amsmath}
\usepackage[utf8]{inputenc}
\usepackage{amssymb}
\usepackage{amsopn}
\usepackage{epsfig}
\usepackage{amsfonts}
\usepackage{latexsym}
\usepackage{graphicx}
\usepackage{calrsfs}
\usepackage{tikz}
\usepackage{enumerate}

\usepackage{dsfont}

 \usetikzlibrary{arrows,automata}
\DeclareMathAlphabet{\pazocal}{OMS}{zplm}{m}{n}
\newtheorem{theorem}{Theorem}[section]
\newtheorem{lemma}[theorem]{Lemma}
\newtheorem{proposition}[theorem]{Proposition}

\theoremstyle{definition}

\theoremstyle{remark}
\newtheorem{remark}[theorem]{Remark}

\numberwithin{equation}{section}

\newcommand{\R}{\ensuremath{\mathbb{R}}}
\newcommand{\N}{\ensuremath{\mathbb{N}}}
\newcommand{\n}{\ensuremath{\vec{n}}}

\newcommand{\us}{\mathbf{U}}

\newcommand{\I}{\mathbf{I}}
\newcommand{\ub}{\mathcal{U}}

\newcommand{\uu}{\pazocal{U}}

\newcommand{\set}[1]{\left\{#1\right\}}
\newcommand{\la}{\lambda}
\newcommand{\ga}{\gamma}
\newcommand{\ep}{\varepsilon}
\newcommand{\f}{\infty}
\newcommand{\de}{\delta}
\newcommand{\om}{\omega}
\newcommand{\al}{\alpha}

\newcommand{\Om}{\Omega}

\newcommand{\lle}{\preccurlyeq}
\newcommand{\lge}{\succcurlyeq}
\begin{document}

\title[Univoque bases of real numbers]{Univoque  bases  of real numbers: simply normal bases, irregular bases and multiple rationals}

\author[Y. Hu]{Yu Hu}
\address[Y. Hu]{College  of Mathematics and Statistics, Chongqing University, Chongqing 401331, People's Republic of China.}
\email{huyu2908@gmail.com}

\author[Y. Huang]{Yan Huang}
\address[Y. Huang]{College  of Mathematics and Statistics, Chongqing University, Chongqing 401331, People's Republic of China.}
\email{yanhuangyh@126.com}

\author[D. Kong]{Derong Kong}
\address[D. Kong]{College  of Mathematics and Statistics, Chongqing University, Chongqing 401331, People's Republic of China.}
\email{derongkong@126.com}

%\date{\today}
\dedicatory{}

%\begin{frontmatter}


\begin{abstract}
Given a positive integer $M$ and a real number $x\in(0,1]$, {we call} $q\in(1,M+1]$ \emph{a univoque simply normal base}  of $x$ if there exists a unique simply normal sequence $(d_i)\in\{0,1,\ldots,M\}^\N$ such that
$x=\sum_{i=1}^\f d_i q^{-i}$. Similarly, a base $q\in(1,M+1]$ is called \emph{a univoque irregular base} of $x$ if there exists a unique sequence $(d_i)\in\{0,1,\ldots, M\}^\N$ such that $x=\sum_{i=1}^\f d_i q^{-i}$ {and the sequence $(d_i)$ has no digit frequency}. Let $\mathcal U_{SN}(x)$ and $\mathcal U_{I_r}(x)$ be the sets of  univoque simply normal bases and   univoque irregular bases of $x$, respectively. In this paper we show that for any $x\in(0,1]$   both $\mathcal U_{SN}(x)$ and $\mathcal U_{I_r}(x)$ have full Hausdorff dimension. Furthermore, given finitely many rationals $x_1, x_2, \ldots, x_n\in(0,1]$ so that each $x_i$ has a  finite {expansion in base $M+1$}, we show that there exists a full {Hausdorff} dimensional set of $q\in(1,M+1]$ such that each $x_i$ has a unique expansion in base $q$.
\end{abstract}
\maketitle

\section{Introduction}\label{sec: Introduction}
Non-integer base expansions were pioneered by R\'enyi \cite{Renyi_1957} and Parry \cite{Parry_1960}. It was extensively studied after the surprising discovery by Erd\H os et al.~\cite{Erdos_Horvath_Joo_1991, Erdos_Joo_1992} that for any $k\in\N\cup\set{\aleph_0}$ there exist $q\in(1,2]$ and $x\in[0, 1/(q-1)]$ such that $x$ has precisely $k$ different $q$-expansions. This phenomenon is completely different from the integer base expansions that each $x$ has a unique expansion, except for countably many $x$ {having} two expansions. In the {literature} of non-integer base expansions there is a great interest in unique   {expansions} due to its close connection with open {dynamical systems} \cite{Erdos_Joo_Komornik_1990, Glendinning_Sidorov_2001} {and kneading theory of unimodal maps or Lorentz maps \cite{Collet-Eckmann-1980, Glendinning-1994, Milnor-Thurston-1988}}.

Given a positive integer $M$ and $q\in(1, M+1]$, a point $x\in[0,  {M}/{(q-1)}]$ is called a \emph{univoque point} in base $q$ if there exists a unique sequence $(d_i)\in\set{0,1,\ldots, M}^\N$ such that
\begin{equation}\label{eq:pi-q}
   x=\sum_{i=1}^\f\frac{d_i}{q^i}=:{((d_i))_q}.
\end{equation}
The infinite sequence {$(d_i)$} is called the (unique) \emph{expansion} of $x$ in base $q$. Throughout the paper we {will fix} the \emph{alphabet} $\set{0,1,\ldots, M}$. Let $\uu_q$ denote the set of all univoque points in base $q$. De Vries and Komornik \cite{DeVries_Komornik_2008, deVries-Komornik-Loreti-2016} studied the topology of $\uu_q$, and showed that $\uu_q$ is closed if and only if $q$ does not belong  to the topological closure of
\[\ub:=\set{p\in(1,M+1]: 1\in\uu_p}.\]
Each base $p\in\ub$ is called a \emph{univoque base} of $1$.
It is known that $\ub$ is a Lebesgue null set of full Hausdorff dimension \cite{Darczy_Katai_1995, DeVries-Komornik-Loreti-2022}. Furthermore, it has a smallest element $\min\ub=q_{KL}$, called the \emph{Komornik-Loreti constant}, which was shown to be transcendental \cite{Allouche_Cosnard_2000, Kong_Li_2015}. Its topological closure $\overline{\ub}$ is a Cantor set {\cite{DeVries-Komornik-Loreti-2022, Komornik_Loreti_2007}}. Some local dimension properties of $\ub$ was studied in \cite{Allaart-Kong-2021}. The set $\ub$ is also related to the bifurcation set of $\alpha$-continued fractions, kneading sequences of unimodal maps, and also  the real slice of the boundary of the Mandelbrot set \cite{Bon-Car-Ste-Giu-2013}.

Motivated by the study of $\ub$,  L\"u, Tan and Wu \cite{Lu_Tan_Wu_2014} initiated the study of  univoque bases of real numbers. Given $x\ge 0$, let
\[
\ub(x):=\set{q\in(1,M+1]: x\in\uu_q}.
\]
Then $\ub=\ub(1)$. Clearly, for $x=0$ we have $\ub(0)=(1,M+1]$, since $0$ {always} has a unique expansion $0^\f=00\ldots$ in any base $q\in(1, M+1]$. When $x\in(0,1]$ {and $M=1$}, L\"u, Tan and Wu showed in \cite{Lu_Tan_Wu_2014} that $\ub(x)$ is a Lebesgue null set but has full Hausdorff dimension. Dajani et al. \cite{Dajani-Komornik-Kong-Li-2018} showed that {for $x \in (0,1]$} the algebraic difference $\ub(x)-\ub(x)$ contains an interval. The third author and his coauthors \cite{Kong-Li-Lv-Wang-Xu-2020} studied the local dimension of $\ub(x)$, and showed that the Hausdorff dimension of $\ub(x)$ is strictly smaller than one {for $x>1$}. Recently, Allaart and the third author \cite{Allaart-Kong-2021} described the smallest element of $\ub(x)$ for {all $x>0$ under the condition $M=1$}.

In this paper we focus on $x\in(0,1]$, and study {those bases of} $\ub(x)$ satisfying some statistical properties.
Given $x\in(0,1]$, a base $q\in(1, M+1]$ is called a \emph{univoque simply normal base} of $x$ if $q\in\ub(x)$ and the unique  expansion $(d_i)\in\set{0,1,\ldots, M}^\N$ of $x$ in base $q$ is simply normal. Let $\ub_{SN}(x)$ be the set of all univoque simply normal bases of $x$. Then for any $q\in\ub_{SN}(x)$   the   unique $q$-expansion $(d_i)$ of $x$  has the same  digit frequency, i.e.,
\[
freq_b((d_i)):=\lim_{n\to\f}\frac{\#\set{1\le i\le n: d_i=b}}{n}=\frac{1}{M+1}\quad\forall ~b\in\set{0,1,\ldots, M}.
\]
Here $\#A$ denotes the cardinality of a set $A$.

On the other hand, a base $q\in(1, M+1]$ is called a \emph{univoque irregular base} of $x$ if $q\in\ub(x)$ and the unique expansion $(d_i)$ of $x$ in base $q$ has no digit frequency. Let $\ub_{I_r}(x)$ be the set of all univoque irregular bases of $x$. {Then for any $q\in \ub_{I_r}(x)$ the unique expansion $(d_i)$ of $x$ in base $q$ satisfies
\begin{equation}\label{eq: inf-fre-ne-sup-fre}
\liminf_{n\to\infty}\frac{\#\{1\leq i\leq n : d_i=b\}}{n}<\limsup_{n\to\infty}\frac{\#\{1\leq i\leq n : d_i=b\}}{n}\quad \forall ~ b\in \{0,1,\ldots M\}.
\end{equation}
{Note} that a sequence $(d_i)$ {satisfying} (\ref{eq: inf-fre-ne-sup-fre}) is  called \emph{essential non-normal} in \cite{Aleberio_Pratsiovytyi_Torbin_2005}.

 It is clear that} $\ub_{I_r}(x)$ and $\ub_{SN}(x)$ are disjoint subsets of $\ub(x)$. %Observe  that $\ub(x)$ is a Lebesgue null set of full Hausdorff dimension for any $x \in (0,1]$.
Our first result states that both {$\ub_{SN}(x)$ and $\ub_{I_r}(x)$} have full Hausdorff dimension for $0<x\le 1$.

\begin{theorem}
  \label{main:simply-normal-irregular}
  For any $x\in(0,1]$ we have
  $
  \dim_H\ub_{SN}(x)=\dim_H\ub_{I_r}(x)=1.
  $
\end{theorem}
Note  that $\ub(x)$ consists of all $q\in(1,M+1]$ such that a given $x$ has a unique $q$-expansion. Then it is natural to ask if we give {finitely many} points $x_1, x_2, \ldots, {x_\ell}>0$, can we find a base $q\in(1, M+1]$ such that each $x_i$ has a unique expansion in base $q$? In general this {$q$ may not exist}, since for example, if $x_1>2$ then the only base such that $x_1$ has a unique expansion is $q=1+M/x_1$ (see \cite[Theorem 1.5]{Kong-Li-Lv-Wang-Xu-2020}). On the other hand, when $x_1, x_2,\ldots, {x_\ell}\in(0,1]$ are all rationals with a finite expansion in base $M+1$, we show that there exists a full {Huasdorff dimensional} set of $q\in(1, M+1]$ such that each $x_i$ has a unique $q$-expansion.

 More precisely, let
\[
D_M:={\left\{\sum_{i=1}^n\frac{d_i}{(M+1)^i}: d_i\in\set{0,1,\ldots, M}\quad\forall~ 1\le i\le n;~ n \in \mathbb{N}\right\}}.
\]
Then $D_M$ is a dense subset of $(0, 1]$. Given finitely many points $x_1, x_2,\ldots, x_\ell\in D_M$, {by constructing}   Cantor subsets of each $\ub(x_i)$ and exploring their {thicknesses}  we show that the intersection  $\bigcap_{i=1}^\ell\ub(x_i)$ has full Hausdorff dimension.
\begin{theorem}
  \label{main:multiple-rationals}
  For any $x_1,\ldots, x_\ell\in D_M$ we have
  \[
  \dim_H\bigcap_{i=1}^\ell\ub(x_i)=1.
  \]
\end{theorem}

The rest of the paper is organized as follows. In Section \ref{sec: preliminaries} we recall some basic properties of univoque bases of real numbers. The proof of Theorem \ref{main:simply-normal-irregular} will be presented in Section \ref{sec:proof of thm1}. To prove Theorem \ref{main:multiple-rationals} we {construct in Section \ref{sec: Cantor-subsets-U(x)} a sequence of Cantor subsets of $\ub(x)$}, and therefore, {by using the thickness} we show {in Theorem \ref{th:algebraic-sum-U(x)}} that the algebraic sum $\ub(x)+\lambda\ub(x)$ contains an interval for any $x\in(0,1]$ and  $\lambda\ne0$. {The proof of  Theorem \ref{main:multiple-rationals}  is in Section \ref{sec:multiple rationals}. In the final section we give some remarks on Theorems \ref{main:simply-normal-irregular} and \ref{main:multiple-rationals}.}

\section{Univoque bases of real numbers}\label{sec: preliminaries}
%{In this section we will introduce some terminology from symbolic dynamics, and recall some properties of unique expansions and univoque bases. This will be used in our proofs of Theorems \ref{main:simply-normal-irregular} and \ref{main:multiple-rationals}. ...}
%
%===========

The study of univoque bases relies  on symbolic dynamics (cf.~\cite{Lind_Marcus_1995}). Given $M\ge 1$, let $\set{0, 1, \ldots, M}^\N$ be the set of all infinite sequences $(d_i)=d_1d_2\ldots$ with each digit $d_i$ from the {{alphabet}} $\set{0,1,\ldots, M}$. By a \emph{word} we mean a finite string of digits over $\set{0,1,\ldots, M}$. Denote by $\set{0,1,\ldots, M}^*$ the set of all finite words including the empty word $\epsilon$. For two words $\mathbf c=c_1\ldots c_m$ and $\mathbf d=d_1\ldots d_n$ we write $\mathbf{cd}=c_1\ldots c_md_1\ldots d_n$ for their concatenation. In particular, for any $k\in\N$ we denote by $\mathbf c^k$ the $k$-fold concatenation of $\mathbf c$ with itself, and by $\mathbf c^\f$ the periodic sequence which is obtained by the infinite concatenation of $\mathbf c$ with itself.

{Throughout} the paper we {will use} lexicographical order $\prec, \lle, \succ$ or $\lge$ between sequences in $\set{0,1,\ldots, M}^\N$. For example, we say  $(i_n)\succ (j_n)$ if $i_1>j_1$, or there exists $n\in\N$ such that $i_1\ldots i_n=j_1\ldots j_n$ and $i_{n+1}>j_{n+1}$. We write $(i_n)\lge (j_n)$ if $(i_n)\succ (j_n)$ or $(i_n)=(j_n)$. Similarly, we say $(i_n)\prec (j_n)$ if $(j_n)\succ (i_n)$, and say $(i_n)\lle (j_n)$ if $(j_n)\lge (i_n)$.
Equipped with {the} metric $\rho$ defined by
\begin{equation}\label{eq:rho}
\rho((i_n), (j_n))=(M+1)^{-\inf\set{n\ge 1: i_n\ne j_n}}
\end{equation}
the symbolic space $\set{0,1,\ldots, M}^\N$ becomes a compact metric space. One can verify that the induced topology by the metric $\rho$ coincides with the order topology on $\set{0,1,\ldots, M}^\N$.

Given $x\in(0,1]$ and $q\in(1, M+1]$, let
\[
\Phi_x(q)=a_1(x,q)a_2(x,q)\ldots\in\set{0,1,\ldots, M}^\N
\]
be the lexicographically largest $q$-expansion of $x$ not ending with $0^\f$, called the \emph{quasi-greedy} $q$-expansion of $x$.
In particular, for   $x=1$  we reserve the notation $\alpha(q)=(\alpha_i(q))$ for the quasi-greedy $q$-expansion of $1$.
The following property  for the quasi-greedy expansion $\Phi_x(q)$   was   proven in \cite[Lemma 2.3 and Lemma 2.5]{DeVries-Komornik-2011}. %For completeness we sketch their main idea and adapt it to our case.
\begin{lemma}
\label{l21}\mbox{}
Let $x\in(0,1]$. The map
$\Phi_x: (1, M+1]\to\set{0,1,\ldots, M}^\N;~ q\mapsto \Phi_x(q)$
 is left continuous under the metric $\rho$, and is strictly increasing   with respect to the lexicographical order. In particular, for $x=1$ the map $q\mapsto \al(q)$ is bijective from $(1, M+1]$ to the set \[\set{(a_i)\in\set{0,1,\ldots, M}^\N: 0^\f\prec a_{n+1}a_{n+2}\ldots \lle a_1a_2\ldots\quad \forall~ n\ge 0 }.\]
\end{lemma}

Given $x\in(0,1]$, recall that $\ub(x)$ consists of all {univoque bases $q\in(1, M+1]$ of $x$}. When $M=1$, the infimum of $\ub(x)$ was {characterized} in \cite{Allaart-Kong-2021}. For a general $M\ge 1$ we still have the following lower bound.
\begin{lemma}
\label{lemma:inf-U(x)}
 For any $x \in (0,1]$ we have
 \[\inf\ub(x)\geq q_G(M)\ge \frac{M}{2}+1,\]
 where \[
  q_G(M):=\left\{\begin{array}
    {lll}
    k+1&\textrm{if}& M=2k,\\
    \frac{k+1+\sqrt{k^2+6k+5}}{2}&\textrm{if}& M=2k+1.
  \end{array}\right.
  \]
\end{lemma}
\begin{proof}
  Note by \cite{Baker_2014} that for $q\in(1, q_G(M)]$ we have $\uu_q=\set{0, M/(q-1)}$. Observe that $q_G(M)\in[\frac{M}{2}+1, M)$. So, {if $x \in (0,1]\cap \uu_q$,} then we must have $q>q_G(M)$. This implies {$\inf \ub(x)\geq q_G(M)$ for any $x \in (0,1]$.}
\end{proof}

For $x\in(0, 1]$ let $$\us(x):=\Phi_x(\ub(x))=\set{\Phi_x(q): q\in\ub(x)}.$$ Then $\Phi_x$ is a bijective map from $\ub(x)$ to $\us(x)$. {Furthermore, the following property of $\Phi_x$ on $\ub(x)$ was shown in \cite[Proposition 3.1 and Proposition 3.3]{Kong-Li-Lv-Wang-Xu-2020}.
\begin{lemma}\label{lem:u(x)-dimension-inequlaity}
  Let $x\in(0,1]$. Then the map $\Phi_x: \ub(x)\to\us(x)$ is locally bi-H\"older continuous under the metric $\rho$ in $\{0,1,\dots, M\}^{\mathbb{N}}$. Furthermore, for any $1<a<b<M+1$ we have
 \[ \frac{\dim_H\Phi_x(\ub(x)\cap(a,b))}{\log b}\le\dim_H(\ub(x)\cap(a,b))\le\frac{\dim_H\Phi_x(\ub(x)\cap(a,b))}{\log a}.\]
\end{lemma}}
Here and throughout the paper we keep using base $M+1$ logarithms.
In view of Lemma \ref{lem:u(x)-dimension-inequlaity}, to study the fractal properties  of $\ub(x)$ and its subsets $\ub_{SN}(x), \ub_{I_r}(x)$ it suffices to study their symbolic analogues
\[\us(x)=\Phi_x(\ub(x)),\quad \us_{SN}(x)=\Phi_x(\ub_{SN}(x))\quad\textrm{and}\quad \us_{I_r}(x)=\Phi_x(\ub_{I_r}(x)).\]
 The following result was essentially obtained in \cite[Section 4]{Lu_Tan_Wu_2014} (see also, \cite[Lemma 4.2]{Kong-Li-Lv-Wang-Xu-2020}).
\begin{lemma}
\label{lem:subset-u-x}
Given $x\in(0,1]$, let $(\ep_i)=\Phi_x(M+1)$ be the quasi-greedy expansion of $x$ in base $M+1$. Then there exist a word $\mathbf w$, a non-negative integer $N$ and a strictly increasing sequence $\set{N_j}_{j=1}^\f\subset\N$ such that
\[
\us_{N_j}(x)\subset\us(x)\quad \textrm{for all }j\ge 1,
\]
where
\[
\us_{N_j}(x):=\set{\ep_1\ldots \ep_{N+N_j} \mathbf w  d_1d_2\ldots:~ d_{n+1}\ldots d_{n+N_j}\notin\set{0^{N_j}, M^{N_j}}\quad\forall~ n\ge 0}.
\]

In particular, if $x\in D_M$, {that is} $(\ep_i)=\Phi_x(M+1)=\ep_1\ldots \ep_m M^\f$ for some $m\ge 1$, {then} we can choose $\mathbf w=\epsilon, N=m$ and $N_j=m+j$.
\end{lemma}

\section{Univoque simply normal bases and univoque irregular bases}\label{sec:proof of thm1}
Let $x\in(0,1]$. Recall that $\ub_{SN}(x)$ consists of all $q\in\ub(x)$ such that $x$ has a unique $q$-expansion which is simply normal.
Furthermore, we recall {from (\ref{eq: inf-fre-ne-sup-fre})} that $\ub_{I_r}(x)$ consists of all $q\in\ub(x)$ such that $x$ has a unique $q$-expansion {with} no digit frequency. Clearly, $\ub_{SN}(x)$ and $\ub_{I_r}(x)$ are disjoint subsets of $\ub(x)$. Note that $\ub(x)$ is a Lebesgue null set of full Hausdorff dimension. In this section we will prove Theorem \ref{main:simply-normal-irregular} that $\dim_H\ub_{SN}(x)=\dim_H\ub_{I_r}(x)=\dim_H\ub(x)=1$ for all $x\in(0,1]$.

\subsection{Univoque simply normal bases}
First we consider the univoque simply normal bases.
\begin{proposition}
  \label{prop:simply-normal}
For any $x\in(0,1]$ we have $\dim_H\ub_{SN}(x)=1$.
\end{proposition}

Our strategy to prove Proposition {\ref{prop:simply-normal}} is to construct a sequence of subsets $\set{\ub_{SN,j}(x)}_{j=1}^\f$ in $\ub_{SN}(x)$ such that $\dim_H\ub_{SN,j}(x)\to 1$ as $j\to\f$.
In view of Lemma \ref{lem:u(x)-dimension-inequlaity} we can do this construction in the symbolic space. For $j\ge 1$ let $\us_{N_j}(x)$ be the subset of $\us(x)$ defined as in Lemma \ref{lem:subset-u-x}. Without loss of generality we assume $N_1>6M$, since otherwise we can delete the first few terms from the sequence $\set{N_j}_{j=1}^\f$. In the following we construct {for each $j\geq 1$} a subset $\us_{SN,j}(x)$ of $\us_{N_j}(x)\cap\us_{SN}(x)$.

Take {$j\geq 1$. Then  $N_j\ge N_1>6M$}. For $k\ge 0$ let
 \begin{equation}\label{eq:mk}
 m_k=2^k(M+1)\lfloor\frac{N_j}{3}\rfloor,
 \end{equation}
where $\lfloor r \rfloor$ denotes the integer part of a real number $r$. {Now for} $k\geq 0$ let $\mathcal N_k$ be the set of all vectors $\n_k:=(n_{k,0}, n_{k,1},\ldots, n_{k,M})$ satisfying
\[
\sum_{b=0}^M n_{k,b}=m_k,\quad\textrm{and}\quad n_{k,b}\in\set{\frac{m_k}{M+1}, \frac{m_k}{M+1}-1}\quad\forall~ 0\le b<M.
\]
{It} is easy to see that $\#\mathcal N_k=2^M$ for all $k\ge 0$. Furthermore, for any $\n_k\in\mathcal N_k$ we have
\begin{equation}
  \label{eq:sn-1}
  \left|\frac{n_{k,b}}{m_k}-\frac{1}{M+1}\right|\le\frac{M}{m_k}\quad\forall~ b\in\set{0,1,\ldots, M}.
\end{equation}
{So $\frac{\n_k}{m_k}$ is a $(M+1)$-dimension probability vector with each element approximately the same.}
Note that $m_k$ and $\mathcal N_k$ both depend on $j$.
In the following  we define the sets ${D_{j,k}}, k\ge 0$ recursively, which will be used to construct our set {$\us_{SN,j}(x)$}.

First we define ${D_{j,0}}$. For a vector $\n_0=(n_{0,0}, n_{0,1},\ldots, n_{0,M})\in\mathcal N_0$ let
\[
D(\n_0):=\set{d_1\ldots d_{m_0}: \xi_b(d_1\ldots d_{m_0})=n_{0,b}\quad\forall~ b\in\set{0,1,\ldots, M}},
\]
where $\xi_b(\mathbf c)$ denotes the number of digit $b$ in the word $\mathbf c$. {Then $D(\n_0)$ consists of all words   {of} length $m_0$ in which each digit $b$ occurs precisely $n_{0,b}$ times.}
 The set ${D_{j,0}}$ is defined by
\[
{D_{j,0}}:=\bigcup_{\n_0\in\mathcal N_0}D(\n_0).
\]

Next suppose ${D_{j,k-1}}$ has been defined for some $k\ge 1$. We define ${D_{j,k}}$ recursively.  {Note {by (\ref{eq:mk})} that $m_k=2m_{k-1}$.} For $\n_k=(n_{k,0}, n_{k,1},\ldots, n_{k,M})\in\mathcal N_k$ let
\[
D(\n_k):=\set{d_1\ldots d_{m_k}\in {D_{j,k-1}}\times {D_{j,k-1}}: \xi_b(d_1\ldots d_{m_k})=n_{k,b}\quad\forall~ b\in\set{0,1,\ldots, M}},
\]
and set
\[
{D_{j,k}}:=\bigcup_{\n_k\in\mathcal N_k}D(\n_k).
\]
Since $\mathcal N_{{k}}$ {consists of} $2^M$ vectors, one can verify that
\begin{equation}
  \label{eq:sn-Dk}
  \begin{split}
\#{D_{j,k}}&=\sum_{\n_k\in\mathcal N_k}\#D(\n_k)\\
&\ge 2^M\min_{\n_k\in\mathcal N_k}\#D(\n_k)\\
&\ge 2^M \left(\min_{\n_0\in\mathcal N_0}\#D(\n_0)\right)^{2^k} =2^{M}\left(\min_{\n_0\in\mathcal N_0}\dbinom {m_0}{\n_0}\right)^{2^k},
\end{split}
\end{equation}
where $\dbinom {m_0}{\n_0}=\dbinom{m_0}{n_{0,0}, n_{0,1},\ldots, n_{0,M}}$ is a multinomial coefficient, and {the second inequality} holds because each block in $D(\n_k)$ belongs to $({D_{j,0}})^{2^k}$.

{Given $x\in(0,1]$, let $(\ep_i)=\Phi_x(M+1)$.} Based on the sets ${D_{j,k}}, k\ge 0$ we define
\[
\us_{SN, j}(x):=\set{\ep_1\ldots \ep_{N+N_j} \mathbf w \mathbf d_0\mathbf d_1\ldots: \mathbf d_k\in {D_{j,k}}\quad\forall~ k\ge 0},
\]
where $N, N_j$ and ${\mathbf w}$ are   defined as in Lemma \ref{lem:subset-u-x}.
\begin{lemma}
  \label{lem:simply-normal-subset}
  Let $x\in(0,1]$. Then for any $j\ge 1$     we have
 $\us_{SN,j}(x)\subset\us_{N_j}(x)\cap\us_{SN}(x).$
\end{lemma}
\begin{proof}
  Note by our construction  that each word $\mathbf d_k\in {D_{j,k}}$ has length $m_k=2^k m_0$ and $\mathbf d_k\in ({D_{j,0}})^{2^k}$. Furthermore, observe that {for each $\mathbf d_0 \in {D_{j,0}}$} the lengths of consecutive zeros and consecutive $M$s in $\mathbf d_0$ are both bounded by $\lfloor\frac{N_j}{3}\rfloor+M$, {which is strictly smaller than the length of $\mathbf d_0$}. So, the lengths of consecutive zeros and consecutive $M$s in each $\mathbf d_k\in {D_{j,k}}$ should be bounded by
  \[
  2\left(\lfloor\frac{N_j}{3}\rfloor+M\right)\le\frac{2N_j}{3}+2M<N_j,
  \]
  where the last inequality holds since   $N_j>6M$. Therefore, by Lemma \ref{lem:subset-u-x} it follows that $\us_{SN,j}(x)\subset\us_{N_j}(x)$.

  To complete the proof we only need to show that each sequence {in} $\us_{SN,j}(x)$ has equal digit frequency. Note that the digit frequency of a sequence is determined by its tail sequences. So, by using (\ref{eq:sn-1}) and that $m_k\to\f$ as $k\to\f$ one can verify that $\us_{SN,j}(x)\subset\us_{SN}(x)$.
\end{proof}

\begin{proof}
  [Proof of Proposition \ref{prop:simply-normal}]
 Let ${\gamma_j}=\max\ub_{N_j}(x)$. {Note by Lemma \ref{lem:subset-u-x}  that {$\Phi_x({\gamma_j})\nearrow\Phi_x(M+1)$ as $j\rightarrow \infty$}. Then by Lemmas \ref{l21} and \ref{lem:u(x)-dimension-inequlaity} it gives that ${\gamma_j}\nearrow M+1$ as $j\to\infty$.} By Lemma \ref{lem:u(x)-dimension-inequlaity}, Lemma \ref{lem:simply-normal-subset} and  \cite[Theorem 2.1]{Feng-Wen-Wu-1997}  it follows that
 \begin{equation}\label{eq:sn-2}
 \begin{split}
   \dim_H\ub_{SN}(x)&\ge\frac{\dim_H\us_{SN}(x)}{\log {\gamma_j}}\ge \frac{\dim_H\us_{SN,j}(x)}{\log {\gamma_j}}\\
    &=\liminf_{n\to\f}\frac{\log\prod_{k=0}^n \# {D_{j,k}}}{\sum_{k=0}^n m_k\log {\gamma_j}}=\liminf_{n\to\f}\frac{\sum_{k=0}^n\log\# {D_{j,k}}}{m_0\sum_{k=0}^n 2^k \log {\gamma_j}}\\
   &\ge \liminf_{n\to\f}\frac{(n+1)\log 2^M+\sum_{k=0}^n 2^k\log\dbinom{m_0}{\n_0^*}}{m_0\sum_{k=0}^n 2^k \log {\gamma_j}}\\
   &=\frac{\log\dbinom{m_0}{\n_0^*}}{m_0\log {\gamma_j}},
 \end{split}
 \end{equation}
 where the last inequality follows by (\ref{eq:sn-Dk}) and $\dbinom{m_0}{\n_0^*}{:=}\min_{\n_0\in\mathcal N_0}\dbinom{m_0}{\n_0}$.
 Note that $\dbinom{m_0}{\n_0^*}=\frac{m_0!}{(n_{0,0}^*! )(n_{0,1}^*!)\ldots(n_{0,M}^*!)}$ {and $m_0=m_0(j)=(M+1)\lfloor\frac{N_j}{3}\rfloor \to\infty$ as $j\to\infty$}. By using {$\sum_{b=0}^M n_{0,b}^*=m_0$ and} the {Stirling's formula that $\log n!=n\log n-n+O(\log n)$ as $n\to\infty$,} it follows that
 \begin{equation}\label{eq:sn-3}
   \begin{split}
     \frac{\log\dbinom{m_0}{\n_0^*}}{m_0\log {\gamma_j}}&= \frac{\log (m_0!)-\sum_{b=0}^M\log(n_{0,b}^*!)}{m_0\log {\gamma_j}}\\
     &=\frac{m_0\log m_0-\sum_{b=0}^M n_{0,b}^*\log n_{0,b}^*+O(\log m_0)}{m_0\log {\gamma_j}}\\
     &=\frac{1}{\log\ga_j}\left(-\sum_{b=0}^M \frac{n_{0,b}^*}{m_0}\log\frac{n_{0,b}^*}{m_0}+O\Big(\frac{\log m_0}{m_0}\Big)\right).
   \end{split}
 \end{equation}
 Observe by (\ref{eq:sn-1}) that for any $b\in\set{0,1,\ldots, M}$,
 \[
 \left|\frac{n_{0,b}^*}{m_0}-\frac{1}{M+1}\right|\le \frac{M}{m_0}=\frac{M}{(M+1)\lfloor\frac{N_j}{3}\rfloor}\to 0\quad\textrm{as }j\to\f.
 \]
Furthermore,
  {${\gamma_j}\to M+1$} as $j\to\f$. So by (\ref{eq:sn-2}) and (\ref{eq:sn-3}) we conclude that
 \[
 \dim_H\ub_{SN}(x)\ge\frac{1}{\log\ga_j}\left(-\sum_{b=0}^M \frac{n_{0,b}^*}{m_0}\log\frac{n_{0,b}^*}{m_0}+O\Big(\frac{\log m_0}{m_0}\Big)\right) \to 1\quad\textrm{as } j\to\f.
 \]
 This completes the proof.
\end{proof}

\subsection{Univoque irregular bases} Now  we consider  the univoque irregular bases.
\begin{proposition}
  \label{prop:irregular-base}
  For any $x\in(0,1]$ we have $\dim_H\ub_{I_r}(x)=1$.
\end{proposition}
{To prove Proposition \ref{prop:irregular-base} we will} construct a sequence of subsets of $\ub_{I_r}(x)$ whose Hausdorff dimension can be arbitrarily close to one. In view of Lemma \ref{lem:u(x)-dimension-inequlaity} it suffices to construct    subsets in $\us_{I_r}(x)=\Phi_x(\ub_{I_r}(x))$.
Recall from Lemma \ref{lem:subset-u-x} that  for any $j\ge 1$,
\[
\us_{N_j}(x)=\set{\ep_1\ldots \ep_{N+N_j}\mathbf w d_1d_2\ldots: d_{n+1}\ldots d_{n+N_j}\notin\set{0^{N_j}, M^{N_j}}\quad \forall ~n\ge 0}
\]
is a subset of $\us(x)$.

 {First we assume $M\ge 2$.  For $k\ge 0$} let ${\Delta_{j,k}}$ be the set of all length $2^k(M+1)N_j(N_j+1)$ words of the form
\begin{equation}\label{eq:Delta-Nj-k}
  c_1\ldots c_{2^k(M+1)N_j^2}\;(0^{N_j-1}1)^{2^k}\;{(1^{N_j})^{2^k} } \cdots {((M-1)^{N_j})^{2^k}}\;(M^{N_j-1}(M-1))^{2^k},
\end{equation}
where $c_i\in\set{0,1,\ldots, M}$ for all $i$, and  $c_{i}\notin\set{0, M}$ if {$i=N_jn$ for some $n \in \mathbb{N}$}. Then each block in ${\Delta_{j,k}}$  has neither  $N_j$ consecutive zeros nor $N_j$ consecutive $M$s. Furthermore,
\begin{equation}\label{eq:number-Delta-Nj-k}
\#{\Delta_{j,k}}=(M+1)^{2^k(M+1)N_j(N_j-1)}(M-1)^{2^k(M+1)N_j}\quad \forall~ k\ge 0.
\end{equation}
{This is because each digit $c_i$ has $M-1$ choices if   the index $i$ is a multiple of $N_j$, and otherwise  $c_i$ has $M+1$ choices.}

Now we define the subset $\us_{I_r,j}(x)$  of $\us_{N_j}(x)$ by
\begin{equation}
  \label{eq:irregular-sequences-j}
  \us_{I_r,j}(x):=\set{\ep_1\ldots \ep_{N+N_j}\mathbf w \mathbf b_0\mathbf b_1\mathbf b_2\cdots: \mathbf b_k\in {\Delta_{j,k}}\quad\forall~ k\ge 0},
\end{equation}
where each ${\Delta_{j,k}}$ is defined in (\ref{eq:Delta-Nj-k}). Since $M\ge 2$, each block $\mathbf b_k$ ends with $M-1\notin\set{0, M}$.
Thus, each sequence $\mathbf b_0\mathbf b_1\cdots$ {contains} neither $N_j$ consecutive zeros nor $N_j$ consecutive $M$s. So, $\us_{I_r,j}(x)$ is indeed a subset of $\us_{N_j}(x)$.

\begin{lemma}
  \label{lem:divergence}
  Let $x\in(0,1]$ and $M\ge 2$. Then for  any $j\ge 1$ we have $\us_{I_r,j}(x)\subset\us_{I_r}(x)$.
\end{lemma}
\begin{proof}
  Note by Lemma \ref{lem:subset-u-x} that $\us_{I_r,j}(x)\subset\us_{N_j}(x)\subset\us(x)$. So it suffices to prove that any {sequence} in $\us_{I_r,j}(x)$ does not have a digit frequency. {Taking a} sequence $\ep_1\ldots \ep_{N+N_j}\mathbf w\mathbf b_0\mathbf b_1\ldots\in \us_{I_r,j}(x)$, we only need to prove that the {tail} sequence
 $
  (d_i)=\mathbf b_0\mathbf b_1\ldots
  $
  has no digit frequency.
  Take a digit $b\in\set{0,1,\ldots, M}$. For $n\in\N$ let
  \[\xi_b(n):=\xi_b(d_1\ldots d_n)=\#\set{1\le i\le n: d_i=b}.\] We will show that the limit of {the sequence  $\{\frac{\xi_b(n)}{n}\}_{n=1}^{\infty}$} does not exist.

  Observe by {(\ref{eq:Delta-Nj-k}) and (\ref{eq:irregular-sequences-j})} that each block $\mathbf b_k$ has length $2^k(M+1)N_j(N_j+1)$ and can be written as
  \[
  \mathbf b_k=\mathbf c_k\, (0^{N_j-1}1)^{2^k}{(1^{N_j})^{2^k}}\cdots {((M-1)^{N_j})^{2^k}}(M^{N_j-1}(M-1))^{2^k},
  \]
  where
  \[\mathbf c_k=c_1 c_2 \ldots c_{2^k(M+1)N_j^2} \qquad\textrm{with } c_i \notin\set{0, M}\textrm{ if } i=N_jn \textrm{ for some }n \in \mathbb{N}.\]
   Let $(\ell_k)$ and $(n_k)$ be two subsequences of $\N$ such that $\ell_k{=\ell_k(b)}$ and $n_k{=n_k(b)}$ are the lengths of blocks $\mathbf b_0\ldots \mathbf b_{k-1}\mathbf c_k (0^{N_j-1}1)^{2^k} {(1^{N_j})^{2^k}}\ldots {((b-1)^{N_j})^{2^k}}$ and $\mathbf b_0\ldots \mathbf b_{k-1}\mathbf c_k (0^{N_j-1}1)^{2^k}{(1^{N_j})^{2^k}}$
   $\ldots {((b-1)^{N_j})^{2^k}}{(b^{N_j})^{2^k}}$, respectively. Then
  \begin{equation}\label{eq:lk-nk}
  \begin{split}
  \ell_k&=2^k(M+1)N_j^2+2^k N_j b+\sum_{i=0}^{k-1}2^i(M+1)N_j(N_j+1),\\
   n_k&=2^k(M+1)N_j^2+2^k N_j (b+1)+\sum_{i=0}^{k-1}2^i(M+1)N_j(N_j+1).
  \end{split}
  \end{equation}
 Furthermore,   let $\theta_b(k)$ be the number of digit $b$ {appearing} in the {block}  $\mathbf c_0\mathbf c_1\ldots \mathbf c_k$. By our definition of $\mathbf c_i$ we must have
 \begin{equation}
   \label{eq:theta-k}
   \begin{split}
   \theta_b(k)&\le\sum_{i=0}^k 2^i(M+1)N_j(N_j-1)={(2^{k+1}-1)}(M+1)N_j(N_j-1)\quad\textrm{if}\quad b\in\set{0, M};\\
   \theta_b(k)&\le \sum_{i=0}^k 2^i(M+1)N_j^2={(2^{k+1}-1)}(M+1)N_j^2\quad\textrm{if}\quad b\in\set{1,2,\ldots, M-1}.
   \end{split}
 \end{equation}

In view of our construction of $(d_i)=\mathbf b_0\mathbf b_1\ldots$, we will finish our proof by considering  the following three cases: (I) $b\in\set{0, M}$; (II) {$b\in\set{1, M-1}$; (III) $b\in\set{2,3,\ldots, M-2}$.}

Case (I). $b\in\set{0, M}$. Then by (\ref{eq:lk-nk}) and the definition of $(d_i)=\mathbf b_0\mathbf b_1\ldots$ it follows that
\begin{equation}\label{eq:tau-lk}
\xi_b(\ell_k)=\theta_b(k)+\sum_{i=0}^{k-1}2^i(N_j-1)=\theta_b(k)+(N_j-1)(2^k-1).
\end{equation}
Suppose on the contrary that $\lim_{n\to\f}\frac{\xi_b(n)}{n}$ exists. Then by (\ref{eq:lk-nk}) and (\ref{eq:tau-lk}) the following limit
\begin{equation}\label{eq:limit-lk}
\begin{split}
  \lim_{k\to\f}\frac{\xi_b(\ell_k)}{\ell_k}&=\lim_{k\to\f}\frac{\theta_b(k)+(N_j-1)(2^k-1)}{2^k(M+1)N_j^2+2^k N_j b+(2^k-1)N_j(M+1)(N_j+1)}\\
  &=\frac{\lim_{k\to\f}\frac{\theta_b(k)}{2^k N_j}+\frac{N_j-1}{N_j}}{(M+1)(2N_j+1)+b}
\end{split}
\end{equation}
exists, which implies that the limit
\[\lim_{k\to\infty}\frac{\theta_b(k)}{2^k N_j}=:A_b\quad \textrm{exists for }b\in\set{0, M}.\]
 Similarly,
\[
\xi_b(n_k)=\theta_b(k)+\sum_{i=0}^{k}2^i(N_j-1)=\theta_b(k)+(N_j-1)(2^{k+1}-1),
\]
and then by (\ref{eq:lk-nk}) and (\ref{eq:limit-lk}) it follows that
\begin{align*}
\frac{A_b+\frac{N_j-1}{N_j}}{(M+1)(2N_j+1)+b}&= \lim_{k\to\f}\frac{\xi_b(\ell_k)}{\ell_k}=\lim_{k\to\f}\frac{\xi_b(n_k)}{n_k}\\
&=\lim_{k\to\f}\frac{\theta_b(k)+(N_j-1)(2^{k+1}-1)}{2^k(M+1)N_j^2+2^k N_j (b+1)+(2^k-1)N_j(M+1)(N_j+1)}\\
&=\frac{A_b+2\frac{N_j-1}{N_j}}{(M+1)(2N_j+1)+b+1},
\end{align*}
which implies that
\[
A_b=\frac{N_j-1}{N_j}\Big[(M+1)(2N_j+1)+b-1\Big]\ge\frac{N_j-1}{N_j}\Big[(M+1)(2N_j+1)-1\Big].
\]This leads to a contradiction, since by (\ref{eq:theta-k}) we have $A_b\le 2(M+1)(N_j-1)$ for $b\in\set{0, M}$.

  Case (II). $b\in\set{1, M-1}$. {First we assume $M\geq 3$}. Then by (\ref{eq:lk-nk}) it follows that
  \[
  \xi_1(\ell_k)=\theta_1(k)+\sum_{i=0}^{k-1}2^i(N_j+1)+2^k,\quad \xi_{M-1}(\ell_k)=\theta_{M-1}(k)+\sum_{i=0}^{k-1}2^i(N_j+1);
  \]
  and
  \[
  \xi_1(n_k)=\theta_1(k)+\sum_{i=0}^k 2^i(N_j+1),\quad \xi_{M-1}(n_k)=\theta_{M-1}(k)+\sum_{i=0}^k 2^i(N_j+1)-2^k.
  \]
  {Suppose the limit $\lim_{n\to\f}\frac{\xi_b(n)}{n}$ exists for $b\in\set{1, M-1}$. Then the limit $A_b:=\lim_{k\to\f}\frac{\theta_b(k)}{2^k N_j}$ also exists.}
  By (\ref{eq:lk-nk}) and the same argument as in Case (I) it follows that
  \[
  \frac{A_1+\frac{N_j+2}{N_j}}{(M+1)(2N_j+1)+1}=\lim_{k\to\f}\frac{\xi_1(\ell_k)}{\ell_k}=\lim_{k\to\f}\frac{\xi_1(n_k)}{n_k}=\frac{A_1+\frac{N_j+2}{N_j}+1}{(M+1)(2N_j+1)+2}
  \]
  and
  \begin{align*}
  \frac{A_{M-1}+\frac{N_j+1}{N_j}}{(M+1)(2N_j+1)+M-1}&=\lim_{k\to\f}\frac{\xi_{M-1}(\ell_k)}{\ell_k} =\lim_{k\to\f}\frac{\xi_{M-1}(n_k)}{n_k}=\frac{A_{M-1}+\frac{N_j+1}{N_j}+1}{(M+1)(2N_j+1)+M}.
  \end{align*}
  This leads to a contradiction with
  (\ref{eq:theta-k}) that $A_b \le 2(M+1)N_j$ for $b\in\set{1, M-1}$. %Hence,  the limit of the sequence  $\{\frac{\tau_n((d_i), b)}{n}\}_{n=1}^{\infty}$ does not exist.

{Next we consider $M=2$. Then $b=1$. By} (\ref{eq:lk-nk}) it follows that
  \[
  \xi_1(\ell_k)=\theta_1(k)+\sum_{i=0}^{k-1}2^i(N_j+2)+2^k \quad  \xi_1(n_k)=\theta_{1}(k)+\sum_{i=0}^k 2^i(N_j+2)-2^k.
  \]
  By (\ref{eq:lk-nk}) and the same argument as above it follows that
  {
  \[
  \frac{A_1+\frac{N_j+3}{N_j}}{(M+1)(2N_j+1)+1}=\lim_{k\to\f}\frac{\xi_1(\ell_k)}{\ell_k}=\lim_{k\to\f}\frac{\xi_1(n_k)}{n_k}=\frac{A_1+\frac{N_j+3}{N_j}+1}{(M+1)(2N_j+1)+2}.
  \]
  Again this leads to a contradiction with  (\ref{eq:theta-k}) that $A_1 \le 2(M+1)N_j$.}
   %$A_1=\lim_{k\to\f}\frac{\theta_1(k)}{2^k N_j}> 2(M+1)N_j$. Hence,  the limit of {the sequence $\{\frac{\tau_n((d_i), b)}{n}\}_{n=1}^{\infty}$} does not exist.

  Case (III). {$b\in\set{2,3,\ldots, M-2}$. Then $M\geq 4$. By} (\ref{eq:lk-nk}) we obtain that
  \begin{align*}
    \xi_b(\ell_k)=\theta_b(k)+\sum_{i=0}^{k-1}2^i N_j,\quad \xi_b(n_k)=\theta_b(k)+\sum_{i=0}^k 2^i N_j.
  \end{align*}
  {Suppose on the contrary that the limit $\lim_{n\to\f}\frac{\xi_b(n)}{n}$ exists for $b\in\set{2,3,\ldots, M-2}$. Then the limit $A_b:=\lim_{k\to\f}\frac{\theta_b(k)}{2^k N_j}$ exists.}
By the same argument as in Case (I) and using (\ref{eq:lk-nk}) we obtain that
\[
\frac{A_b+1}{(M+1)(2N_j+1)+b}=\lim_{k\to\f}\frac{\xi_b(\ell_k)}{\ell_k}=\lim_{k\to\f}\frac{\xi_b(n_k)}{n_k}=\frac{A_b+2}{(M+1)(2N_j+1)+b+1},
\]
{which leads to a contradiction, since by (\ref{eq:theta-k}) we have $A_b\le 2(M+1)N_j$ for $b\in\set{2,\ldots,  M-2}$.

Therefore, by Cases (I)--(III) we conclude that the frequency of digit $b$ in $(d_i)=\mathbf b_0\mathbf b_1\mathbf b_2\ldots$ does not exist for any $b\in\set{0,1,\ldots, M}$. This completes the proof.}
%we can also prove that the limit of {the sequence $\{\frac{\tau_n((d_i), 1)}{n}\}_{n=1}^\infty$} does not exist.
  \end{proof}

  \begin{proof}
    [Proof of Proposition \ref{prop:irregular-base}]First we consider $M\ge 2$.
    Let $x\in(0,1]$ and let $(\ep_i)=\Phi_x(M+1)$. Suppose $\us_{N_j}(x), j=1,2,\ldots$ are the subsets {of $\us(x)$} defined as in Lemma \ref{lem:subset-u-x}. {For $j\ge 1$ let  $\gamma_j$ be the largest element of $\ub_{N_j}(x):=\Phi_x^{-1}(\us_{N_j}(x))$. Then $\Phi_x(\gamma_j)=\ep_1\ldots\ep_{N+N_j}\mathbf w (M^{N_j-1}(M-1))^\f\nearrow\Phi_x(M+1)$ as $j\to\f$. So by Lemmas \ref{l21} and \ref{lem:u(x)-dimension-inequlaity} it follows that $\gamma_j\nearrow M+1$ as $j\to\f$.}

    Note by Lemma \ref{lem:divergence} that $\us_{I_r,j}(x)\subset\us_{I_r}(x)$. So, by Lemma \ref{lem:u(x)-dimension-inequlaity} it follows that
    \begin{equation}\label{eq:irregular-dim}
      \dim_H\ub_{I_r}(x)\ge \dim_H\ub_{I_r,j}(x) \ge\frac{\dim_H\us_{I_r,j}(x)}{\log {\gamma_j}}
    \end{equation}
    for all $j\ge 1$. Note that $(\set{0,1,\ldots, M}^\N, \rho)$ is a compact metric space, where $\rho$ is defined in (\ref{eq:rho}). Then by {(\ref{eq:Delta-Nj-k}), (\ref{eq:number-Delta-Nj-k}) and \cite[Theorem 2.1]{Feng-Wen-Wu-1997} }it follows that
    \begin{align*}
      \dim_H\us_{I_r,j}(x)&= \liminf_{n\to\f}\frac{\log\prod_{k=0}^n\#{\Delta_{j,k}}}{\sum_{k=0}^n|\mathbf b_k|\log(M+1)}\\
      &=\liminf_{n\to\f}\frac{\sum_{k=0}^n 2^k(M+1)N_j\Big[(N_j-1)\log(M+1)+\log(M-1)\Big]}{\sum_{k=0}^n 2^k(M+1)N_j(N_j+1)\log(M+1)}\\
      &=\frac{N_j-1}{N_j+1}+\frac{\log(M-1)}{(N_j+1)\log(M+1)}.
    \end{align*}
    Since ${\gamma_j}\to  M+1$ and $N_j\to\f$ as $j\to\f$, by (\ref{eq:irregular-dim}) this implies that
    \[
    \dim_H\ub_{I_r}(x)\ge \frac{N_j-1}{(N_j+1)\log {\gamma_j}}+\frac{\log(M-1)}{(N_j+1)\log(M+1)\log {\gamma_j}}\to1 \quad \text{as }j\to\infty.
    \]
    Here we emphasize that the logarithm is in base $M+1$.

%\medskip

Now we consider  $M=1$. The proof is similar.   We modify the definition of $\us_{I_r,j}(x)$ as
\[
\widetilde\us_{I_r,j}(x)=\set{\ep_1\ldots \ep_{{N+N_j}}\mathbf w \mathbf b_0\mathbf b_1\ldots: \mathbf b_k\in{\widetilde\Delta_{j,k}}\quad\forall ~k\ge 0},
\]
where each ${\widetilde\Delta_{j,k}}$ consists of all length $2^{k+1}N_j(N_j+1)+2$ blocks of the form
\[
 c_1\ldots c_{2^{k+1}N_j^2}\;(0^{N_j-1}1)^{2^k}(01^{N_j-1})^{2^k}\,01
\]
with {each $c_i\in\set{0,1}$, and $c_{i-1}c_i=01$ if $i=N_j n$ for some $n\in\N$.} Then each sequence $\mathbf b_0\mathbf b_1\ldots\in\prod_{k=0}^\f{\widetilde\Delta_{j,k}}$  contains neither $N_j$ consecutive zeros nor $N_j$ consecutive {ones}. So, $\widetilde\us_{I_r,j}(x)\subset\us_{N_j}(x)$. By the same argument as in the proof of Lemma  \ref{lem:divergence} one can verify that $\widetilde\us_{I_r,j}(x)\subset\us_{I_r}(x)$. Hence, by similar argument as above we can prove that $\dim_H\ub_{I_r}(x)=1$ for $M=1$.
\end{proof}
\begin{proof}
  [Proof of Theorem \ref{main:simply-normal-irregular}] The theorem follows by Propositions \ref{prop:simply-normal} and \ref{prop:irregular-base}.
\end{proof}

\section{Cantor subsets of $\ub(x)$ and thickness}\label{sec: Cantor-subsets-U(x)}
 In this section we will show that {$\ub(x)+\la\ub(x)$ contains an interval for any $x\in(0,1]$ and $\la\ne 0$}, which generalizes the main result of \cite{Dajani-Komornik-Kong-Li-2018} where they proved this result only  {for $M=1$.} %In fact, we prove the following more general result.

\begin{theorem}
  \label{th:algebraic-sum-U(x)}
If $f: \R^2\to \R$ is $C^1$   such that the partial derivatives are not vanishing in $(M+1-\de, M+1]^2$ for some $\de>0$, then for any $x\in(0,1]$ the set
  \[
\ub_f(x):=\set{f(p, q): p, q\in\ub(x)}
  \]
  contains an interval.
\end{theorem}
\begin{remark}\mbox{}
  \begin{enumerate}[{\rm(i)}]
  \item If $f(x,y)=x+\lambda y$ for some $\lambda\ne 0$, then $f$ is $C^1$ with its partial derivatives not vanishing. So, by Theorem \ref{th:algebraic-sum-U(x)} it follows that ${\ub}(x)+\lambda {\ub}(x)$ contains an interval {for any $x\in(0,1]$}.
    \item For   possible extension of this theorem we refer to a recent paper   \cite{Jiang-2022} and the references therein.
  \end{enumerate}
\end{remark}

\subsection{Thickness of a Cantor set in $\R$}
The thickness of a Cantor set in $\R$ was introduced by Newhouse \cite{Newhouse-1979}, and it has {been applied} in dynamical systems and number theory (cf.~\cite{Astels-2000}). Let $E\subset\R$ be a Cantor set with the convex hull $E_0$. Then the complement $E_0\setminus E=\bigcup_{n=1}^\f O_n$ is the union of countably many disjoint open intervals. The sequence $\mathcal O=(O_1, O_2,O_3,\ldots)$ is called a \emph{derived sequence} of $E$. If the lengths of these open intervals are in {a} non-increasing order, i.e., $|O_i|\ge|O_{i+1}|$ for all $i\ge 1$, then we call the sequence $\mathcal O$ an \emph{ordered sequence}. Let $E_n:=E_0\setminus\bigcup_{k=1}^n O_k$. Then for any $n\ge 1$ the open interval $O_n$ must belong to {a unique} connected component ${C}$ of $E_{n-1}$. {In this case,} ${C}\setminus O_n$ is the union of two {disjoint} closed intervals $L_{\mathcal O}(O_n)$ and $R_{\mathcal O}(O_n)$. Hence, the thickness of $E$ with respect to the derived sequence $\mathcal O$ is defined by
\[
\tau_{\mathcal O}(E):=\inf_{n\ge 1}\min\set{\frac{|L_{\mathcal O}(O_n)|}{|O_n|}, \frac{|R_{\mathcal O}(O_n)|}{|O_n|}}{;}
\]
{and} the \emph{thickness} of $E$ is {then} defined by
\begin{equation}
  \label{eq:def-thickness}
  \tau(E):=\sup\set{\tau_{\varrho(\mathcal O)}(E): \varrho(\mathcal O)\textrm{ is a permutation of }\mathcal O}.
\end{equation}
Note by \cite{Astels-2000} that the supremum in (\ref{eq:def-thickness}) is attainable. Indeed, for any ordered sequence $\mathcal O$ we have $\tau(E)=\tau_{\mathcal O}(E)$.

The following result for the relationship between the thickness of a Cantor set and its Hausdorff dimensioin was given by Newhouse \cite{Newhouse-1979} (see also, \cite{Palis_Takens_1993}).
\begin{lemma}
  \label{lem:thickness-Hausdorff-dim}
  Let $E\subset \R$ be a Cantor set. Then
  \[
  \dim_H E\ge\frac{\log 2}{\log(2+1/\tau(E))}.
  \]
\end{lemma}
From Lemma \ref{lem:thickness-Hausdorff-dim} it follows that if the thickness of a Cantor set $E$ is very large, then its Hausdorff dimension is close to $1$.
The next result{,} which can be derived from \cite{Hunt-Kan-Yorke-1993}{,} {describes} how the thickness can be used to study the intersection of two Cantor sets. %We call two Cantor sets in $\R$ \emph{interleaved} if neither set lies in a gap of the other.
\begin{lemma}\label{lem:intersection-thickness}
Let $E$ and $F$ be two  Cantor sets in $\R$ having the same maximum point $\xi$. If $\xi$ is an accumulation point of $E\cap F$, and their {thicknesses} $\tau(E)\ge t$ and $\tau(F)\ge t$ for some large $t>0$, then there exists a Cantor subset $K\subset E\cap F$ such that $\max K=\xi$ and $\tau(K)\ge C \sqrt{t}$ for some  $C>0$.
\end{lemma}

 The following result on the image of two Cantor sets $E$ and $F$ can be deduced from \cite{McDonald-Taylor-2021} and \cite{Simon-Taylor-2020} (see also, \cite{Jiang-2022}).
 \begin{lemma}
   \label{lem:image-Cantor-sets}
   Let $E$ and $F$ be two Cantor sets in $\R$ with $\tau(E)\tau(F)>1$.  If $f:\R^2\to\R$ is a $C^1$ function with non-vanishing partial derivatives, then the set $\set{f(x, y): x\in E, y\in F}$ contains an interval.
 \end{lemma}

\subsection{Proof of Theorem \ref{th:algebraic-sum-U(x)}}
Fix $x\in(0,1]$ {let $(\ep_i)=\Phi_x(M+1)$. Recall from Lemma \ref{lem:subset-u-x} that $\set{\us_{N_j}(x)}_{j\ge 1}$ is a sequence of subsets in $\us(x)$.}   Then for any $j\ge 1$, each sequence of $\us_{N_j}(x)$ ends neither with $N_j$ consecutive zeros nor $N_j$ consecutive $M$s. Set
\begin{equation}\label{eq:u-N-j-x-q}
\ub_{N_j}(x)=\{q\in (1,M+1]:\Phi_x(q)\in \us_{N_j}(x)\}.
\end{equation}
Then by Lemma \ref{l21}  it follows that $\Phi_x:\ub_{N_j}(x)\to\us_{N_j}(x)$ is an increasing homeomorphism. Observe that each $\us_{N_j}(x)$ is a Cantor set with respect to the metric $\rho$ defined in (\ref{eq:rho}). This implies that $\ub_{N_j}(x)$ is a Cantor subset of $(1,M+1]$, and so it can be geometrically constructed by successively removing a sequence of open intervals from a closed interval.

To describe the geometrical construction of $\ub_{N_j}(x)$ we first define a sequence of symbolic intervals. For $j\ge 1$ let
\[
\Om_j^*(x):=\bigcup_{n=0}^\f\Om_j^n(x),
\]
where
\[
\Om_j^n(x):=\set{\ep_1\ldots\ep_{N+N_j}\mathbf w d_1d_2\ldots d_n: d_{i+1}\ldots d_{i+N_j}\notin\set{0^{N_j}, M^{N_j}}\quad\forall ~0\le i\le n-N_j}.
\]
Here $N, N_j$ and $\mathbf w$ are defined as in Lemma \ref{lem:subset-u-x}.
%Then $\Om_j^*(x)$ consists of all words constructed by concatenating  $\vs_j$ with admissible words of the subshift   of finite type which has forbidden blocks $0^{N_j}$ and $M^{N_j}$ (cf.~\cite{Lind_Marcus_1995}).
Now for each $\om\in\Om_j^*(x)$ we denote by
$\I_\om=[(a_i), (b_i)]$   the symbolic interval {which contains all sequences in $\us_{N_j}(x)$ beginning with $\om$. Here  $(a_i)$ and $(b_i)$ are the lexicographically smallest and largest sequences in $\us_{N_j}(x)$ beginning with $\om$, respectively.} %So, for $\om=\vs_j\in\Om_0(x)$ we have \[\I_{\vs_j}=[\vs_j(0^{N_j-1}1)^\f, \vs_j(M^{N_j-1}(M-1))^\f].\] For a general word $\om\in\Om_j^*(x)$ we have the following.

\begin{lemma}\label{lem:symbolic interval}
Let $x\in(0,1]$ and $j\ge 1$. Take $\omega =\ep_1\ldots\ep_{N+N_j}\mathbf w \;\mathbf d\in \Omega^*_j(x)$.
\mbox{}
\begin{enumerate}[{\rm(i)}]
\item If $\mathbf d$ ends with neither $0$ nor $M$, then
\[
\mathbf{I}_{\omega}=[\omega(0^{N_j-1}1)^\infty,\ \omega(M^{N_j-1}(M-1))^{\infty}].
\]
\item  If $\mathbf d$ ends with $0^k$ {for some} $k\in \{ 1,2,\ldots,N_j-1\}$, then
\[
\mathbf{I}_{\omega}=[\omega0^{N_j-1-k}(10^{N_j-1})^{\infty},\ \omega(M^{N_j-1}(M-1))^{\infty}].
\]
\item If $\mathbf d$ ends with $M^k$ {for some} $k\in \{ 1,2,\ldots,N_j-1\}$, then
\[
\mathbf{I}_{\omega}=[\omega(0^{N_j-1}1)^\infty,\ \omega M^{N_j-1-k}((M-1)M^{N_j-1})^\infty].
\]
\end{enumerate}
\end{lemma}

\begin{figure}[h!]
\begin{center}
\begin{tikzpicture}[
    scale=12,
    axis/.style={very thick, ->},
    important line/.style={thick},
    dashed line/.style={dashed, thin},
    pile/.style={thick, ->, >=stealth', shorten <=2pt, shorten
    >=2pt},
    every node/.style={color=black}
    ]
    % axis

    % Lines
   % \draw[important line]  (0,0.85)--( 1.0, 0.85);
%    \node[] at(0.5, 0.9){$I_\om$};
    \draw[important line] (0, 0)--(1.1, 0);
     \draw[important line] (0, -0.1) --(0.18, -0.1);  \node[above,scale=0.8pt] at(0.09, -0.15){$I_{\om 0}$}; \node[above,scale=0.8pt] at(0.23, -0.1){$G_{\om 0}$};  \node[above,scale=0.8pt] at(0.34, -0.15){$I_{\om 1}$};\draw[important line] (0.28, -0.1) --(0.4, -0.1);
     %\draw[important line] (0, -0.2) --(0.08, -0.2);  \node[above,scale=0.7pt] at(0.1, -0.2){$V_{\ell,3}$}; \draw[important line] (0.12, -0.2) --(0.18, -0.2);
%     \draw[important line] (0.28, -0.2) --(0.33, -0.2);  \node[above,scale=0.7pt] at(0.345, -0.2){$V_{\ell,2}$}; \draw[important line] (0.36, -0.2) --(0.4, -0.2);

     \node[above,scale=1pt]at(0.55,0.02){$I_\om$};

\draw[dashed line] (0.43, -0.1) --(0.52, -0.1);

   % \draw[important line] (0.55, 0)node[above,scale=1pt]{$\al_{\ell+1}$}--(0.8, 0)node[above,scale=1pt]{$\beta_{\ell+1}$};
    \draw[important line] (0.55, -0.1) --(0.66, -0.1);  \node[above,scale=0.8pt] at(0.6, -0.15){$I_{\om d}$}; \node[above,scale=0.8pt] at(0.68, -0.1){$G_{\om d}$};\node[above,scale=0.8pt] at(0.75, -0.15){$I_{\om(d+1)}$}; \draw[important line] (0.7, -0.1) --(0.8, -0.1);
    % \draw[important line] (0.55, -0.2) --(0.595, -0.2);  \node[above,scale=0.7pt] at(0.61, -0.2){$V_{\ell+1,3}$}; \draw[important line] (0.625, -0.2) --(0.66, -0.2);
%     \draw[important line] (0.7, -0.2) --(0.745, -0.2);  \node[above,scale=0.7pt] at(0.755, -0.2){$V_{\ell+1,2}$}; \draw[important line] (0.765, -0.2) --(0.8, -0.2);

     % \node[above,scale=1pt]at(0.675,-0.35){$F_{\ell+1}(x)$};

\draw[dashed line] (0.82, -0.1) --(0.88, -0.1);

    %\draw[important line] (0.9, 0)node[above,scale=1pt]{$\al_{\ell+2}$}--(1.1, 0)node[above,scale=1pt]{$\beta_{\ell+2}$};
    \draw[important line] (0.9, -0.1) --(0.99, -0.1);  \node[above,scale=0.8pt] at(0.945, -0.15){$I_{\om(M-1)}$};\node[above,scale=0.8pt] at(1.01, -0.1){$G_{\om(M-1)}$};\node[above,scale=0.8pt] at(1.065, -0.15){$I_{\om M}$}; \draw[important line] (1.03, -0.1) --(1.1, -0.1);
     %\draw[important line] (0.9, -0.2) --(0.94, -0.2);  \node[above,scale=0.7pt] at(0.95, -0.2){$V_{\ell+2,3}$}; \draw[important line] (0.96, -0.2) --(0.99, -0.2);
%     \draw[important line] (1.03, -0.2) --(1.06, -0.2);  \node[above,scale=0.7pt] at(1.07, -0.2){$V_{\ell+2,2}$}; \draw[important line] (1.08, -0.2) --(1.1, -0.2);
 %\node[above,scale=1pt]at(1,-0.35){$F_{\ell+2}(x)$};
\end{tikzpicture}
\end{center}
\caption{The $(n+1)$-level gaps $G_{\om d}, d=0,1,\ldots, M-1$ associated to the $n$-level basic interval $I_\om$ in the construction of $\ub_{N_j}(x)$.}\label{fig:1}
\end{figure}
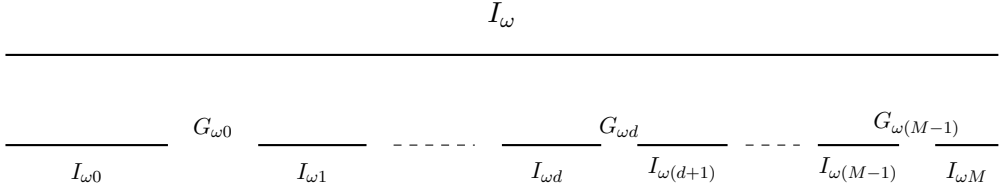

Now we describe the geometrical construction of $\ub_{N_j}(x)$ in terms of the symbolic intervals $\set{\I_\om: \om\in\Om_j^*(x)}$. For a symbolic interval $\I_\om=[(a_i), (b_i)]$ with $\om\in\Om_j^*(x)$ we define the associated interval $I_\om=[p, q]$ by
\[
 \Phi_x(p)=(a_i)\quad\textrm{and}\quad \Phi_x(q)=(b_i).
\]
Since $\Phi_x$ is an increasing homeomorphism from $\ub_{N_j}(x)$ to $\us_{N_j}(x)$, it follows that the convex hull of $\ub_{N_j}(x)$ is $I_{\ep_1\ldots\ep_{N+N_j}\mathbf w}$. Furthermore, for any $n\ge 0$ and any $\om\in\Om_j^n(x)$ the intervals $I_{\om d}, \om d\in\Om_j^{n+1}(x)$ are pairwise disjoint subintervals of $I_\om$. It turns out that the set $\set{I_\om: \om\in\Om_j^*(x)}$ of basic intervals  has a {tree} structure. Therefore,
\begin{equation}\label{eq:geometic-u-j(x)}
\ub_{N_j}(x)=\bigcap_{n=0}^\f\bigcup_{\om\in\Om_j^n(x)}I_\om.
\end{equation}
Each closed interval $I_\om$ with $\om\in\Om_j^n(x)$ is called an \emph{$n$-level basic interval}. We emphasize that the endpoints of each $n$-level basic interval belong to $\ub_{N_j}(x)$. Now for a $n$-level basic interval $I_\om$  we define the \emph{$(n+1)$-level gaps} associated to $I_\om$ as follows: suppose $I_{\om d}$ and $I_{\om(d+1)}$ are two consecutive {$(n+1)$-level} basic intervals, then the gap between them, denoted by $G_{\om d}$, is a $(n+1)$-level gap (see Figure \ref{fig:1}). By Lemma \ref{lem:symbolic interval} it follows that the number of $(n+1)$-level gaps associated to $I_\om$ is either $M$ or $M-1$, and the later case refers to items (ii) and (iii) in Lemma \ref{lem:symbolic interval}.

Based on the geometrical construction of $\ub_{N_j}(x)$ in (\ref{eq:geometic-u-j(x)}) it is convenient to  define its thickness {according to the basic intervals $I_\om$ and gaps $G_\om$ for $\om\in\Om_j^*(x)$. Let}
\begin{equation}\label{eq:thickness-new}
\tau_*(\ub_{N_j}(x))=\inf_{n\ge 0}\min_{\om\in\Om_j^n(x)}\set{\frac{|I_{\om }|}{|G_{\om}|},
\frac{|I_{\om^+}|}{|G_{\om}|}},
\end{equation}
where for {$\om\in \Om_j^n(x)$, if $\om^+ \notin\Om_j^n(x)$} we set $\frac{|I_{\om^+}|}{|G_{\om}|}=+\f$. Here for a word $\mathbf c=c_1\ldots c_k$ we set $\mathbf c^+:=c_1\ldots c_{k-1}(c_k+1)$.
 By (\ref{eq:def-thickness}) and (\ref{eq:thickness-new}) it follows that $\tau(\ub_{N_j}(x))\ge \tau_*(\ub_{N_j}(x))$.

\begin{proposition}
  \label{prop:algebraic-sum-interval}
  For any $x\in(0,1]$ we have $\tau_*(\ub_{N_j}(x))\to\f$ as $j\to\f$.
\end{proposition}

\begin{proof}
  Let $x\in(0,1]$. Then $(\ep_i)=\Phi_x(M+1)\succ 0^\f$, and there exists $\ell\in\N$ such that $\ep_1\ldots \ep_\ell\lge 0^{\ell-1}1$. So, we can take $j\in\N$ large enough such that
  \begin{equation}\label{eq:ell}
     \varepsilon_1\ldots \varepsilon_{N+N_j} \mathbf w 0^\infty \succ 0^{\ell-1}10^\f .
  \end{equation}
   {Then each sequence in $\us_{N_j}(x)$ is lexicographically larger than $0^{\ell-1}1 0^\f$.} {For $n\ge 0$ let} $\om=\ep_1\ldots\ep_{N+N_j}\mathbf w d_1\ldots d_n\in\Om_j^n(x)$. {For brevity we write $n_j:=|\omega|=N+N_j+|\mathbf w|+n$ for its length.} Suppose $\om d, \om(d+1)\in\Om_j^{n+1}(x)$. Write
  \[
  I_{\om d}=[q_1, q_2],\quad I_{\om(d+1)}=[q_3, q_4].
  \]
  Then the gap $G_{\om d}=(q_2, q_3)$. By Lemma \ref{lem:symbolic interval} it follows that
  \begin{equation}
    \label{eq:28-1}
    \begin{split}
      &\om d(0^{N_j-1}1)^\f\lle \Phi_x(q_1)\lle \om d(10^{N_j-1})^\f,\\
      &\Phi_x(q_2)=\om d(M^{N_j-1}(M-1))^\f{,}\\
      &\Phi_x(q_3)=\om(d+1)(0^{N_j-1}1)^\f,\\
      &\om(d+1)((M-1)M^{N_j-1})^\f\lle\Phi_x(q_4)\lle \om (d+1)(M^{N_j-1}(M-1))^\f.
    \end{split}
  \end{equation}
  Then
  \[
  (\om d(M^{N_j-1}(M-1))^\f)_{q_2}=x=(\om(d+1)(0^{N_j-1}1)^\f)_{q_3},
  \]
  which together with (\ref{eq:ell}) implies
  \begin{equation}\label{eq:28-2}
  \begin{split}
   &{(0^{{n_j}}1(0^{N_j-1}1)^\f)_{q_3}-(0^{{n_j}+1}(M^{N_j-1}(M-1))^\f)_{q_2}}\\
    =&~(\om d 0^\f)_{q_2}-(\om d 0^\f)_{q_3} \ge\frac{1}{q_{{2}}^\ell}-\frac{1}{q_{{3}}^\ell}
    \geq  {\frac{q_3-q_2}{q_2q_3^\ell}} \ge \frac{q_3-q_2}{q_3^{\ell+1}}.
  \end{split}
  \end{equation}
  On the other hand,
  \begin{equation}\label{eq:28-3}
  \begin{split}
    &({0^{{n_j}}}1(0^{N_j-1}1)^\f)_{q_3}-({0^{{n_j+1}}}(M^{N_j-1}(M-1))^\f)_{q_2}\\
    \le &~({0^{{n_j}}}1(0^{N_j-1}1)^\f)_{q_3}-({0^{{n_j+1}}}(M^{N_j-1}(M-1))^\f)_{q_3}\\
    \le &~({0^{{n_j+1}}}(M^{N_j-1}(M+1))^\f)_{q_3}-({0^{{n_j+1}}}(M^{N_j-1}(M-1))^\f)_{q_3}\\
    =&~\frac{2}{q_3^{{n_j+N_j+1}}(1-1/q_3^{N_j})}<\frac{4}{q_3^{{n_j+N_j+1}}},
  \end{split}
  \end{equation}
  where the second inequality follows by $({0^{{n_j+1}}}M^\f)_{q_3}\geq({0^{{n_j}}}10^\f)_{q_3}$, and the last inequality holds since $1/q_3^{N_j}<1/2$ for large $j$. Hence, by (\ref{eq:28-2}) and (\ref{eq:28-3}) we {obtain} an upper bound on the length of $G_{\om d}$:
  \begin{equation}
    \label{eq:gap-upper}
    |G_{\om d}|=q_3-q_2{<}\frac{4}{q_3^{{n_j+N_j-\ell}}}.
  \end{equation}

 In the following we consider the lower bounds on the lengths of $I_{\om d}$ and $I_{\om(d+1)}$. {To do this we need} the following inequalities.\\
 {\bf Claim:} for all sufficiently large {$j$}  we have
 \begin{equation}\label{eq:28-**}
 \begin{split}
 ({0^{{n_j+1}}}(10^{N_j-1})^\f)_{q_1}&\le ({0^{{n_j+1}}}M(M-1)(M^{N_j-3}(M-1)M^2)^\f)_{q_2},\\
 ({0^{{n_j+1}}}01(0^{N_j-3}10^2)^\f)_{q_3}&\le ({0^{{n_j+1}}}((M-1)M^{N_j-1})^\f)_{q_4}.
 \end{split}
 \end{equation}

 {Since} the proofs of the two inequalities in (\ref{eq:28-**}) are similar, we only prove the first inequality. Note by (\ref{eq:28-1}) that
 \[
 (\om d(0^{N_j-1}1)^\f)_{q_1}\le x=(\om d(M^{N_j-1}(M-1))^\f)_{q_2}.
 \]
 This together with (\ref{eq:ell}) implies that
\begin{equation}
   \begin{split}
     ({0^{{n_j+1}}}M^\f)_{q_2}&\ge ({0^{{n_j+1}}}(M^{N_j-1}(M-1))^\f)_{q_2}-({0^{{n_j+1}}}(0^{N_j-1}1)^\f)_{q_1}\\
     &{\geq} (\om d0^\f)_{q_1}-(\om d 0^\f)_{q_2} \ge \frac{1}{q_1^\ell}-\frac{1}{q_2^\ell}\ge\frac{q_2-q_1}{q_1q_2^\ell}.
   \end{split}
 \end{equation}
 Whence,
 \begin{equation}\label{eq:inequality-1}
 \frac{q_2-q_1}{q_1}\le\frac{M}{q_2^{{n_j+1}-\ell}(q_2-1)}\le\frac{M}{q_G^{{n_j}-\ell}},
 \end{equation}
 where the last inequality holds since $q_i\ge q_G=q_G(M)$ by Lemma \ref{lemma:inf-U(x)}. Therefore, the first inequality of (\ref{eq:28-**}) can be deduced as follows:
 \begin{align*}
   (0^{{{n_j+1}}}(10^{N_j-1})^\f)_{q_1}&=\left(1+\frac{q_2-q_1}{q_1}\right)^{{{n_j+1}}}\frac{((10^{N_j-1})^\f)_{q_1}}{q_2^{{{n_j+1}}}}\\
   &\le \left(1+\frac{M}{{q_G^{{n_j}-\ell}}}\right)^{{{n_j+1}}}\frac{((10^{N_j-1})^\f)_{q_G}}{q_2^{{{n_j+1}}}}\\
   &\le \frac{(M(M-1)(M^{N_j-3}(M-1)M^2)^\f)_{M+1}}{q_2^{{{n_j+1}}}}\\
   &\le (0^{{{n_j+1}}}M(M-1)(M^{N_j-3}(M-1)M^2)^\f)_{q_2},
 \end{align*}
  where the first inequality follows by (\ref{eq:inequality-1}), and the second inequality holds for all sufficiently large {$j$} since
  \[
  \lim_{j\to\f}\left(1+\frac{M}{{q_G^{{n_j}-\ell}}}\right)^{{{n_j+1}}}=1<  \lim_{j\to\f}\frac{(M(M-1)(M^{N_j-3}(M-1)M^2)^\f)_{M+1}}{((10^{N_j-1})^\f)_{q_G}}.
  \]
This proves the claim.

\medskip
Now by (\ref{eq:28-1}) we have
\[
(\om d(10^{N_j-1})^\f)_{q_1}\ge x=(\om d(M^{N_j-1}(M-1))^\f)_{q_2},
\]
which implies that
\begin{equation}
  \label{eq:28-4}
  \begin{split}
    &(0^{{{n_j+1}}}(M^{N_j-1}(M-1))^\f)_{q_2}-(0^{{{n_j+1}}}(10^{N_j-1})^\f)_{q_1}\\
    \le& ~(\om d0^\f)_{q_1}-(\om d 0^\f)_{q_2}{\le} (M^\f)_{q_1}-(M^\f)_{q_2}\\
    =& ~{\frac{M}{(q_1-1)(q_2-1)}(q_2-q_1)}\leq \frac{M}{(q_G-1)^2}(q_2-q_1).
  \end{split}
\end{equation}
 On the other hand, by the first inequality of (\ref{eq:28-**}) it follows that
\begin{equation}
\label{eq:28-5}
\begin{split}
&(0^{{{n_j+1}}}(M^{N_j-1}(M-1))^\f)_{q_2}-(0^{{{n_j+1}}}(10^{N_j-1})^\f)_{q_1}\\
\ge&~(0^{{{n_j+1}}}(M^{N_j-1}(M-1))^\f)_{q_2}-(0^{{{n_j+1}}}M(M-1)(M^{N_j-3}(M-1)M^2)^\f)_{q_2}\\
= &~(0^{{n_j}+2}10^\f)_{q_2}=\frac{1}{q_2^{{n_j}+3}}.
\end{split}\end{equation}
So, by (\ref{eq:28-4}) and (\ref{eq:28-5}) we obtain a lower bound on the length of $I_{\om d}$:
\begin{equation}
  \label{eq:inter-lower-d}
  |I_{\om d}|=q_2-q_1\ge \frac{(q_G-1)^2}{M q_2^{{n_j+3}}}> \frac{(q_G-1)^2}{M q_3^{{n_j+3}}}.
\end{equation}

Similarly, by (\ref{eq:28-1}) we also have
\[
(\om(d+1)(0^{N_j-1}1)^\f)_{q_3}=x\ge(\om(d+1)((M-1)M^{N_j-1})^\f)_{q_4},
\]
which implies that
\begin{equation}
  \label{eq:28-7}
  \begin{split}
    &(0^{{{n_j+1}}}((M-1)M^{N_j-1})^\f)_{q_4}-(0^{{{n_j+1}}}(0^{N_j-1}1)^\f)_{q_3}\\
    \le &~(\om(d+1)0^\f)_{q_3}-(\om(d+1)0^\f)_{q_4}{\le} \frac{M}{(q_G-1)^2}(q_4-q_3).
  \end{split}
\end{equation}
On the other hand, by the second inequality of (\ref{eq:28-**}) it follows that
\begin{equation}
  \label{eq:28-8}
  \begin{split}
  &(0^{{{n_j+1}}}((M-1)M^{N_j-1})^\f)_{q_4}-(0^{{{n_j+1}}}(0^{N_j-1}1)^\f)_{q_3}\\
  \ge&~(0^{{{n_j+1}}}01(0^{N_j-3}10^2)^\f)_{q_3}-(0^{{{n_j+1}}}(0^{N_j-1}1)^\f)_{q_3}{=}\frac{1}{q_3^{{n_j}+3}}.
  \end{split}
\end{equation}
  So, by (\ref{eq:28-7}) and (\ref{eq:28-8}) we obtain a lower bound on the length of $I_{\om(d+1)}$:
  \begin{equation}
    \label{eq:inter-lower-d+1}
    |I_{\om(d+1)}|=q_4-q_3\ge\frac{(q_G-1)^2}{M q_3^{{n_j}+3}}.
  \end{equation}

Hence, by (\ref{eq:gap-upper}), (\ref{eq:inter-lower-d}) and (\ref{eq:inter-lower-d+1}) we conclude that
\[
\min\set{\frac{|I_{\om d}|}{|G_{\om d}|}, \frac{|I_{\om(d+1)}|}{|G_{\om d}|}}{\ge}\frac{(q_G-1)^2}{4M}q_3^{N_j-\ell-3} \to \;\f\quad \textrm{as } j\to\f.
\]
 Since $\om$ was taken from $\Om_j^n(x)$ arbitrarily, this completes the proof by (\ref{eq:thickness-new}).
\end{proof}
\begin{proof}
  [Proof of Theorem \ref{th:algebraic-sum-U(x)}]
 Take $x\in(0,1]$. Note by Proposition \ref{prop:algebraic-sum-interval} that
  \[
 \tau(\ub_{N_j}(x))\ge \tau_*(\ub_{N_j}(x))\to \f \quad \textrm{as }  j\to\f.
 \]
 Let $f:\R^2\to\R$ be a $C^1$ function with partial derivatives not vanishing on $(M+1-\de, M+1]^2$ for some $\de>0$. Observe that $\ub_{N_j}(x)\subset(M+1-\de, M+1]$ for all large $j$. Then by Lemma \ref{lem:image-Cantor-sets} it follows that $\ub_f(x)=\set{f(p,q): p, q\in\ub(x)}$ contains an interval.
\end{proof}

\section{Univoque bases of multiple rationals}\label{sec:multiple rationals}
{In this section we will prove Theorem \ref{main:multiple-rationals}. Recall that $D_M$ consists of all rationals in $[0,1]$ with a finite expansion in base $M+1$. Given $x_1, x_2,\ldots, x_{{\ell}}\in D_M$, we will show that the intersection $\bigcap_{i=1}^{{\ell}}\ub(x_i)$ has full Hausdorff dimension.}

Take $x\in D_M$. Then $\Phi_x(M+1)=\ep_1\ldots \ep_m M^\f$ for some $\ep_m<M$ {with $m\geq1$.  {Let $j\ge m$.} Then $2^j>m$, and by Lemma \ref{lem:subset-u-x} it follows that
\[
{\us_{2^j}}(x)\subset\us(x),
\]
where
\[
{\us_{2^j}}(x)=\set{\ep_1\ldots\ep_{m} M^{{2^j}} d_1d_2\ldots: d_{n+1}\ldots d_{{n+{2^j}}}\notin\set{0^{{{{2^j}}}}, M^{{{{2^j}}}}}\quad\forall ~n\ge 0}.
\]
Accordingly, ${\ub_{2^j}}(x)=\Phi_x^{-1}({\us_{2^j}}(x))\subset \ub(x)$. Let $[\al_j, \beta_j]$ be the convex hull of ${\ub_{2^j}}(x)$. In the following we show that these subintervals $[\al_j,\beta_j], j\ge m$ are pairwise disjoint and {converge} to $\set{M+1}$ under the {Hausdorff metric} (see Figure \ref{fig:2}).
   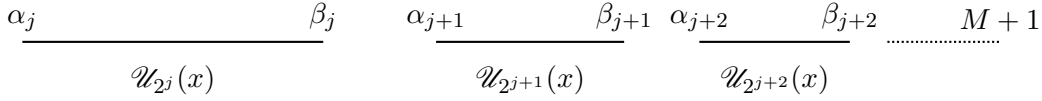
\begin{figure}[h!]
\begin{center}
\begin{tikzpicture}[
    scale=10,
    axis/.style={very thick, ->},
    important line/.style={thick},
    dashed line/.style={dashed, thin},
    pile/.style={thick, ->, >=stealth', shorten <=2pt, shorten
    >=2pt},
    every node/.style={color=black}
    ]
    % axis

    % Lines
   % \draw[important line]  (0,0.85)--( 1.0, 0.85);
%    \node[] at(0.5, 0.9){$I_\om$};
    \draw[important line] (0, 0)node[above,scale=1pt]{$\al_j$}--(0.4, 0)node[above,scale=1pt]{$\beta_j$};
     \node[] at(0.2, -0.05){${\ub_{2^j}}(x)$};

    \draw[important line] (0.55, 0)node[above,scale=1pt]{$\al_{j+1}$}--(0.8, 0)node[above,scale=1pt]{$\beta_{j+1}$};

\node[] at(0.675, -0.05){${\ub_{2^{j+1}}}(x)$};

    \draw[important line] (0.9, 0)node[above,scale=1pt]{$\al_{j+2}$}--(1.1, 0)node[above,scale=1pt]{$\beta_{j+2}$};
    \node[] at(1, -0.05){${\ub_{2^{j+2}}}(x)$};

  \draw[important line, densely dotted] (1.15,0)--(1.3,0)node[above,scale=1pt]{$M+1$};

\end{tikzpicture}
\end{center}
\caption{The geometry  of the convex hulls $[\al_j, \beta_j]=\textrm{conv}({\ub_{2^j}}(x)),~ j\ge m$.}\label{fig:2}
\end{figure}
\begin{lemma}
  \label{lem:alpha-j-beta-j}
  Let $x\in D_M$ and $[\al_j, \beta_j]=\textrm{conv}({\ub_{2^j}}(x))$ with $j\ge m$. Then
  \[
 \al_j< \beta_j<\al_{j+1}\quad\forall~ j\ge m,\quad\textrm{and}\quad \al_j\nearrow M+1\quad\textrm{as }j\to\f.
  \]
\end{lemma}
\begin{proof}
 Take $j\ge m$. By the definition of $\ub_{2^j}(x)$ it follows that
  \begin{equation}
    \label{eq:alpha-beta-j}
    \begin{split}
      \Phi_x(\al_j)&=\ep_1\ldots \ep_{m}M^{{2^j}}(0^{{{2^j-1}}}1)^\f,\\
      \Phi_x(\beta_j)&=\ep_1\ldots \ep_{m}M^{{2^j}}(M^{{{2^j-1}}}(M-1))^\f,\\
      \Phi_x(\al_{j+1})&=\ep_1\ldots \ep_{m}M^{2^{j+1}}(0^{2^{j+1}-1}1)^\f,\\
      \Phi_x(M+1)&=\ep_1\ldots \ep_{m}M^\f.
    \end{split}
  \end{equation}
  Then $\Phi_x(\al_j)\prec \Phi_x(\beta_j)\prec \Phi_x(\al_{j+1})$, and thus $\al_j<\beta_j<\al_{j+1}$ by Lemma \ref{l21}. Furthermore, $\Phi_x(\al_j)$ increasingly converges to $\ep_1\ldots \ep_m M^\f=\Phi_x(M+1)$ under the metric $\rho$. {By Lemma \ref{l21} and Lemma \ref{lem:u(x)-dimension-inequlaity}} it follows that $\al_j\nearrow M+1$ as $j\to\f$.
\end{proof}

Note that each ${\ub_{2^j}}(x)$ is a Cantor set, i.e., a nonempty compact set has neither interior nor isolated points. So, by Lemma \ref{lem:alpha-j-beta-j} it follows that for any {$k\ge m$},
\[
E_k(x):=\bigcup_{j=k}^\f{\ub_{2^j}}(x)\cup\set{M+1}
\]
is also a Cantor set.

\begin{proposition}
\label{prop:algebraic-intersect-interval}
  For any $x\in D_M$ we have $\tau(E_k(x))\to\f$ as $k\to\f$.
\end{proposition}
\begin{proof}
  Since $E_k(x)=\bigcup_{j=k}^\f{\ub_{2^j}}(x)\cup\set{M+1}$ is a Cantor set, by (\ref{eq:def-thickness}) it suffices to prove that for each $k\in\N$ there exists a derived sequence $\mathcal O_k$ of $E_k(x)$ such that $\tau_{\mathcal O_k}(E_k(x))\to\f$ as $k\to\f$. Observe that
  \[
  \textrm{conv}(E_k(x))\setminus E_k(x)=\bigcup_{j=k}^\f(\beta_j, \al_{j+1})\cup\bigcup_{j=k}^\f\bigcup_{\ell=1}^\f O_\ell^{(j)},
  \]
  where $\mathcal O^{(j)}=(O_1^{(j)}, O_2^{(j)}, O_3^{(j)},\ldots)$ is an ordered sequence of $\ub_{2^j}(x)$. By Proposition \ref{prop:algebraic-sum-interval}  it follows that
  \begin{equation}
    \label{eq:thickness-infinity-uj}
    \tau(\ub_{2^j}(x))=\tau_{\mathcal O^{(j)}}(\ub_{2^j}(x))\to\f\quad\textrm{as }j\to\f.
  \end{equation}

  Now we define a derived sequence $\mathcal O_k$ of $E_k(x)$ which consists of open intervals ordered in the following way:
  \begin{align*}
    &(\beta_k,\al_{k+1}),\quad O_1^{(k)};\\
    &(\beta_{k+1}, \al_{k+2}), \quad O_2^{(k)}, \quad O_1^{(k+1)};\\
    &(\beta_{k+2}, \al_{k+3}),\quad O_3^{(k)}, \quad O_2^{(k+1)},\quad  O_1^{(k+2)};\\
    &\cdots;\\
    &(\beta_{k+n-1}, \beta_{k+n}),\quad O_n^{(k)}, \quad O_{n-1}^{(k+1)},\quad\ldots, \quad O_1^{(k+n-1)};\\
    &\cdots.
  \end{align*}
 % Note by (\ref{eq:def-thickness}) that to prove $\lim_{k\to\f}\tau(E_k(x))=\f$ it suffices to prove
%  \[\lim_{k\to\f} \tau_{\mathcal O_k} (E_k(x))=\f.\]
  So, by our definition of the derived sequence $\mathcal O_k$ and  (\ref{eq:thickness-infinity-uj}) it suffices to prove that
  \begin{equation}
    \label{eq:thickness-al-j-beta-j}
    \min\set{\frac{\beta_j-\al_j}{\al_{j+1}-\beta_j}, \frac{M+1-\al_{j+1}}{\al_{j+1}-\beta_j}}\to\f\quad\textrm{as }j\to\f.
  \end{equation}

  First we give an upper bound of $\al_{j+1}-\beta_j$. Note by (\ref{eq:alpha-beta-j}) that
 \[
(\ep_1\ldots \ep_m M^{2^j}(M^{2^{j}-1}(M-1))^\f)_{\beta_j}=x=(\ep_1\ldots \ep_m M^{2^{j+1}}(0^{2^{j+1}-1}1)^\f)_{\al_{j+1}}.
 \]
 Then
 \begin{align*}
 &(0^{m+2^{j+1}-1}1(0^{2^{j+1}-1}1)^\f)_{\al_{j+1}}-(0^{m+2^{j+1}}(M^{2^j-1}(M-1))^\f)_{\beta_j}\\
 =~&(\ep_1\ldots \ep_m M^{2^{j+1}-1}(M-1)0^\f)_{\beta_j}-(\ep_1\ldots\ep_m M^{2^{j+1}-1}(M-1)0^\f)_{\al_{j+1}}\\
 \ge~&\frac{M}{\beta_j^{m+1}}-\frac{M}{\al_{j+1}^{m+1}}
 \ge\frac{M}{\al_{j+1}^{m+2}}(\al_{j+1}-\beta_j).
 \end{align*}
 This implies that
 \begin{align*}
   \al_{j+1}-\beta_j&\le \frac{\al_{j+1}^{m+2}}{M}\left[(0^{m+2^{j+1}-1}1(0^{2^{j+1}-1}1)^\f)_{\al_{j+1}}-(0^{m+2^{j+1}}(M^{2^j-1}(M-1))^\f)_{\beta_j}\right]\\
  %&\le\al_{j+1}^{m+2}\left[(0^{m+2^{j+1}}(M^{2^{j+1}-1}(M+1))^\f)_{\al_{j+1}}-(0^{m+2^{j+1}}(M^{2^j-1}(M-1))^\f)_{\al_{j+1}}\right]\\
   &\le\al_{j+1}^{m+2}\left[(0^{m+2^{j+1}}(M^{2^j-1}(M+1))^\f)_{\al_{j+1}}-(0^{m+2^{j+1}}(M^{2^j-1}(M-1))^\f)_{\al_{j+1}}\right]\\
   &=\al_{j+1}^{m+2}\frac{((20^{2^j-1})^\f)_{\al_{j+1}}}{\al_{j+1}^{m+2^{j+1}+2^j{-1}}}.
 \end{align*}
 Thus, by Lemma \ref{lemma:inf-U(x)} we obtain an upper bound of $\al_{j+1}-\beta_j$:
 \begin{equation}
   \label{eq:upper-bound}
   \al_{j+1}-\beta_j\le \frac{((20^{2^j-1})^\f)_{\al_{j+1}}}{\al_{j+1}^{2^{j+1}+2^j-3}}\le\frac{C_0}{\al_{j+1}^{2^{j+1}+2^j}},
 \end{equation}
 where $C_0:=\frac{2(M+1)^{3}}{q_G-1}$.

 {Next we consider the  lower bounds  of $\beta_j-\al_j$ and $M+1-\alpha_{j+1}$, respectively. These are based on the following two inequalities.\\
 {\bf Claim:} for all sufficiently large $j$  we have
 \begin{equation}\label{eq:sec-5-claim}
 \begin{split}
 (0^{m+2^j}(0^{2^j-1}1)^\f)_{\alpha_j}&\le (0^{m+2^j}0(M^{2^j-2}(M-1)M)^\f)_{\beta_j},\\
 (0^{m+2^{j+1}}M(0^{2^{j+1}-2}10)^\f)_{\alpha_{j+1}}&\le (0^{m+2^{j+1}}M^\f)_{M+1}.
 \end{split}
 \end{equation}

Since the proof of the first inequality in (\ref{eq:sec-5-claim}) is similar to the proof of (\ref{eq:28-**}), we only prove the second inequality. Note by (\ref{eq:alpha-beta-j}) that
$$(\varepsilon_1\ldots \varepsilon_m M^\infty)_{M+1}=x=(\varepsilon_1\ldots \varepsilon_m M^{2^{j+1}}(0^{2^{j+1}-1}1)^\infty)_{\alpha_{j+1}}.$$ Then
\begin{align*}
(0^{m+2^{j+1}}M^\infty)_{M+1}&\geq (0^{m+2^{j+1}}M^\infty)_{M+1}-(0^{m+2^{j+1}}(0^{2^{j+1}-1}1)^\infty)_{\alpha_{j+1}}\\
&=(\varepsilon_1\ldots \varepsilon_m M^{2^{j+1}}0^\infty)_{\alpha_{j+1}}-(\varepsilon_1\ldots \varepsilon_m M^{2^{j+1}}0^\infty)_{M+1}\\
&\geq \frac{1}{\alpha_{j+1}^{m+1}}-\frac{1}{(M+1)^{m+1}}\geq \frac{M+1-\alpha_{j+1}}{\alpha_{j+1}(M+1)^{m+1}}.
\end{align*}
Whence,
\begin{equation}\label{eq:m+1-al-j+1}
\frac{M+1-\alpha_{j+1}}{\alpha_{j+1}}\leq \frac{(M+1)^{m+1}}{(M+1)^{m+2^{j+1}}}=\frac{1}{(M+1)^{2^{j+1}-1}}.
\end{equation}
Therefore, the second inequality in (\ref{eq:sec-5-claim}) can be deduced as follows:
\begin{align*}
(0^{m+2^{j+1}}M(0^{2^{j+1}-2}10)^\infty)_{\alpha_{j+1}}&= \left(1+\frac{M+1-\alpha_{j+1}}{\alpha_{j+1}}\right)^{m+2^{j+1}}\frac{(M(0^{2^{j+1}-2}10)^\infty)_{\alpha_{j+1}}}{(M+1)^{m+2^{j+1}}}\\
&\leq \left(1+\frac{1}{(M+1)^{m+2^{j+1}}}\right)^{m+2^{j+1}}\frac{(M(0^{2^{j+1}-2}10)^\infty)_{\alpha_{j+1}}}{(M+1)^{m+2^{j+1}}}\\
&\leq \frac{1}{(M+1)^{m+2^{j+1}}}=(0^{m+2^{j+1}}M^\infty)_{M+1},
\end{align*}
where the first inequality follows by (\ref{eq:m+1-al-j+1}), and the second inequality holds for all sufficiently large $j$ since
$$\lim_{j\to\infty}\left(1+\frac{1}{(M+1)^{2^{j+1}-1}}\right)^{m+2^{j+1}}=1
<\lim_{j\to\infty}\frac{1}{(M(0^{2^{j+1}-2}10)^\infty)_{\alpha_{j+1}}}.$$
This proves the claim.
}

\medskip

 Note by (\ref{eq:alpha-beta-j}) that
 \[
 (\ep_1\ldots \ep_m M^{2^j}(M^{2^j-1}(M-1))^\f)_{\beta_j}=x=(\ep_1\ldots\ep_m M^{2^j}(0^{2^j-1}1)^\f)_{\al_j}.
 \]
 Then
 \begin{align*}
&(0^{m+2^j}(M^{2^j-1}(M-1))^{\infty})_{\beta_j}-(0^{m+2^j}(0^{2^j-1}1)^{\infty})_{\alpha_j}\\
=~&(\ep_1 \cdots \ep_{m}M^{2^j}0^{\infty})_{\alpha_j}-(\ep_1 \cdots \ep_{m}M^{2^j}0^{\infty})_{\beta_j}\\
\leq~&(M^{\infty})_{\alpha_j}-(M^{\infty})_{\beta_j}=\frac{M}{(\al_j-1)(\beta_j-1)}(\beta_j-\al_j)
\leq\frac{M}{(q_G-1)^2}  \left(\beta_j-\alpha_j\right).
\end{align*}
This implies that
\begin{align*}
\beta_j-\alpha_j &\geq \frac{(q_G-1)^2}{M} \left[(0^{m+2^j}(M^{2^j-1}(M-1))^{\infty})_{\beta_j}-(0^{m+2^j}(0^{2^j-1}1)^{\infty})_{\alpha_j}\right] \\
&\geq \frac{(q_G-1)^2}{M} \left[(0^{m+2^j}(M^{2^j-1}(M-1))^{\infty})_{\beta_j}-(0^{m+2^j}0(M^{2^j-2}(M-1)M)^{\infty})_{\beta_j}\right]
\end{align*}
for sufficiently large $j$,
where the second inequality {follows by the first inequality in (\ref{eq:sec-5-claim})}.
Therefore,
\begin{equation}
\label{eq:lower-bound-1}
\beta_j-\alpha_j \geq  \frac{(q_G-1)^2 }{\beta_j^{m+2^j+1}}\ge\frac{C_1}{\beta_j^{2^j}}>\frac{C_1}{\al_{j+1}^{2^j}},
\end{equation}
where $C_1:=\frac{(q_G-1)^2}{(M+1)^{m+1}}$.

  Now we turn to a lower bound of $M+1-\al_{j+1}$. Note by (\ref{eq:alpha-beta-j}) that
  \[
  (\ep_1\ldots \ep_m M^{2^{j+1}}(0^{2^{j+1}-1}1)^\f)_{\al_{j+1}}=x=(\ep_1\ldots \ep_m M^\f)_{M+1}.
  \]
  Then
  \begin{align*}
&(0^{m+2^{j+1}}M^{\infty})_{M+1}-(0^{m+2^{j+1}}(0^{2^{j+1}-1}1)^{\infty})_{\alpha_{j+1}}\\
=~&(\ep_1 \cdots \ep_{m}M^{2^{j+1}}0^{\infty})_{\alpha_{j+1}}-(\ep_1 \cdots \ep_{m}M^{2^{j+1}}0^{\infty})_{M+1} \\
\leq~& (M^{\infty})_{\alpha_{j+1}}-(M^{\infty})_{M+1}
{\leq \frac{1}{ q_G-1 }   \left(M+1-\alpha_{j+1}\right)}.
\end{align*}
This implies that
\begin{align*}
M+1-\al_{j+1}&\ge (q_G-1)\left[(0^{m+2^{j+1}}M^{\infty})_{M+1}-(0^{m+2^{j+1}}(0^{2^{j+1}-1}1)^{\infty})_{\alpha_{j+1}}\right]\\
&\geq (q_G-1)\left[(0^{m+2^{j+1}}M(0^{2^{j+1}-2}10)^{\infty})_{\alpha_{j+1}}-(0^{m+2^{j+1}}(0^{2^{j+1}-1}1)^{\infty})_{\alpha_{j+1}}\right]
\end{align*}
for $j$ sufficiently large, where the second inequality follows by {the second inequality in (\ref{eq:sec-5-claim})}.
So,
\begin{equation}
\label{eq:lower-bound-2}
M+1-\alpha_{j+1} \geq   \frac{M(q_G-1)}{\alpha_{j+1}^{m+2^{j+1}+1}}\ge\frac{C_2}{\alpha_{j+1}^{2^{j+1}}},
\end{equation}
where $C_2:=\frac{M(q_G-1)}{(M+1)^{m+1}}$.

  Hence, by (\ref{eq:upper-bound}), (\ref{eq:lower-bound-1}) and (\ref{eq:lower-bound-2}) it follows that
  \begin{align*}
    \min\set{\frac{\beta_j-\al_j}{\al_{j+1}-\beta_j}, \frac{M+1-\al_{j+1}}{\al_{j+1}-\beta_j}}&\ge\min\set{\frac{C_1}{C_0}\al_{j+1}^{2^{j+1}}, \frac{C_2}{C_0}\al_{j+1}^{2^j}}\to\f
  \end{align*}
  as $j\to\f$. This proves (\ref{eq:thickness-al-j-beta-j}) and then completes the proof.
\end{proof}

\begin{proof}
  [Proof of Theorem \ref{main:multiple-rationals}]
  Let $x_1, x_2, \ldots, x_\ell\in D_M$. {Then there exists $m\in\N$ such that for each $x_i$ the sequence $\Phi_{x_i}(M+1)=d_1d_2\ldots$ satisfies $d_{m+1}d_{m+2}\ldots=M^\f$.} By our construction each set $\ub(x_i)\cup\set{M+1}$ contains a sequence of Cantor subsets $E_k(x_i)=\bigcup_{j=k}^\f{\ub_{2^j}}(x_i)\cup\set{M+1}, {k\ge m}$. By Proposition {\ref{prop:algebraic-intersect-interval}} the thickness $\tau(E_k(x_i))\to \f$ as $k\to\f$. Furthermore, each $E_k(x_i)$ has a maximum value $M+1$. So, by Lemma \ref{lem:intersection-thickness} {and Lemma \ref{lem:alpha-j-beta-j}} it follows that for any {$k\ge m$} and any $i_1, i_2\in\set{1,2,\ldots, \ell}$ the intersection $E_k(x_{i_1})\cap E_k(x_{i_2})$ contains a Cantor subset $E_k(x_{i_1}, x_{i_2})$ such that
  \[
  \max E_k(x_{i_1}, x_{i_2})=M+1,\quad\textrm{and}\quad \tau(E_k(x_{i_1}, x_{i_2}))\to\f\quad\textrm{as }k\to\f.
    \]

    Proceeding this argument for all $x_1, x_2,\ldots x_\ell\in D_M$ we obtain that for any {$k\ge m$} the intersection $\bigcap_{i=1}^\ell E_k(x_i)$ contains a Cantor subset $E_k(x_1,\ldots, x_\ell)$ satisfying
    \[
    \max E_k(x_1, x_2,\ldots, x_\ell)=M+1,\quad\textrm{and}\quad \tau(E_k(x_1,x_2,\ldots, x_\ell))\to\f\quad\textrm{as }k\to\f.
    \]
    Hence, by Lemma \ref{lem:thickness-Hausdorff-dim} we conclude that
    \begin{align*}
      \dim_H\bigcap_{i=1}^\ell \ub(x_i)&\ge\dim_H\bigcap_{i=1}^\ell {E_k(x_i)}\ge\dim_H E_k(x_1,\ldots, x_\ell)\\
      &\ge\frac{\log 2}{\log\left(2+\frac{1}{\tau(E_k(x_1,\ldots, x_\ell))}\right)}\quad\to 1\quad\textrm{as }k\to\f.
    \end{align*}
 This completes the proof.
\end{proof}

\section{Final remarks}
{In the way of proving Theorem \ref{main:simply-normal-irregular} we can also obtain that
\[\sup\ub_{SN}(x)=\sup\ub_{I_r}(x)=M+1\quad \textrm{for any }x\in(0,1].\] On the other hand, for $M=1$ and $x=1$ the smallest element of $\ub(1)$ is $q_{KL}\approx 1.78723$, and $\Phi_1(q_{KL})$ is the shift of the Thue-Morse sequence (cf.~\cite{Komornik-Loreti-1998}). Note that the Thue-Morse sequence is simply normal (cf.~\cite{Allouche-Shallit-1998}).}
      So, $\min\ub_{SN}(1)=q_{KL}$ for $M=1$. {However, when $M>1$ the smallest element of $\ub(1)$ is not univoque simply normal (cf.~\cite{Komornik_Loreti_2002}). Then it is natural to ask what  $\inf\ub_{NS}(1)$ is for $M>1$. In general, for $x\in(0,1)$ and $M\ge 1$ it is interesting to determine $\inf\ub_{SN}(x)$. Also, it might be interesting to investigate $\inf\ub_{I_r}(x)$ for a general $x\in(0,1]$ and $M\ge 1$.}

 {In Theorem \ref{main:multiple-rationals} we show that the intersection $\bigcap_{i=1}^\ell\ub(x_i)$ has full Hausdorff dimension for any  $x_1,\ldots, x_\ell\in D_M$.  However, if some $x_i$ does not belong to $D_M$,  then our method does not work. So it is worth exploring  whether Theorem \ref{main:multiple-rationals} holds for any finitely many $x_1, x_2,\ldots, x_\ell\in(0,1]$.}

{\section*{Acknowledgements}
The first author was supported by CYS22075. The third author was supported by NSFC No.~11971079.

%\bibliographystyle{abbrv}
%\bibliography{Fractal-Expansions}

\end{document}